\documentclass[a4paper,12pt]{amsart}

\usepackage{float}
\usepackage{euscript,eufrak,verbatim}
\usepackage{graphicx}
\usepackage[usenames]{color}
\usepackage[colorlinks,linkcolor=red,anchorcolor=blue,citecolor=blue]{hyperref}
\usepackage{amsmath}
\usepackage{amsthm}
\usepackage[all]{xy}
\usepackage{graphicx} 
\usepackage{amssymb} 

\usepackage[all]{xy}

\usepackage{mathrsfs}
\usepackage{amscd}

\makeatletter
\numberwithin{equation}{section}

\theoremstyle{plain}
\newtheorem{main theorem}{Main Theorem}
\newtheorem{theorem}{Theorem}[section]
\newtheorem{lemma}[theorem]{Lemma}

\newtheorem{corollary}[theorem]{Corollary}
\newtheorem{proposition}[theorem]{Proposition}
\newtheorem{claim}[theorem]{Claim}

\theoremstyle{definition}

\newtheorem{remark}[theorem]{Remark}
\newtheorem{example}[theorem]{Example}

\newtheorem{problem}[theorem]{Problem}

\newtheorem{question}[theorem]{Question}
\newcommand{\diam}{\mathrm{diam}}

\newcommand{\h}{h_{\mathrm{top}}}

\makeatother

\begin{document}

\title[New approach to weighted topological entropy and pressure]{New approach to weighted topological entropy and pressure}

\author[Masaki Tsukamoto]{Masaki Tsukamoto}



\thanks{}

\maketitle

\begin{abstract}
Motivated by fractal geometry of self-affine carpets and sponges, Feng--Huang (2016)
introduced weighted topological entropy and pressure for factor maps between dynamical systems, 
and proved variational principles for them.
We introduce a new approach to this theory.
Our new definitions of weighted topological entropy and pressure are very different from 
the original definitions of Feng--Huang.
The equivalence of the two definitions seems highly nontrivial.
Their equivalence can be seen as a generalization of the dimension formula for the Bedford--McMullen carpet
in purely topological terms.
\end{abstract}

\section{Introduction} \label{section: introduction}

\subsection{Weighted topological entropy and pressure}  \label{subsection: weighted topological entropy and pressure}

The purpose of this paper is to introduce a new approach to 
\textit{weighted topological entropy and pressure} introduced by 
Feng--Huang \cite{Feng--Huang}.
In this subsection we describe their original theory.
We explain our new approach in the next subsection. 

We first quickly review the classical theory of entropy and pressure of dynamical systems.
See the book of Walters \cite{Walters book} for the details.
A pair $(X, T)$ is called a \textbf{dynamical system} if $X$ is a compact metrizable space and 
$T:X\to X$ is a continuous map.
We denote its \textbf{topological entropy} by $h_{\mathrm{top}}(X,T)$.
This is a topological invariant of dynamical systems, which counts the number of bits per iterate for describing the orbits of 
$(X,T)$.

One of the most basic theorems about topological entropy is \textbf{variational principle}.
We define $\mathscr{M}^T(X)$ as the set of invariant Borel probability measures on $X$.
For each measure $\mu \in \mathscr{M}^T(X)$, we denote its \textbf{Kolomogorov--Sinai entropy} by $h_\mu(T)$.
Then the variational principle states that \cite{Goodwyn, Dinaburg, Goodman}
\begin{equation} \label{eq: standard variational principle}
  h_{\mathrm{top}}(X, T) = \sup_{\mu\in \mathscr{M}^T(X)} h_\mu(T). 
\end{equation}  

This theory can be generalized to \textit{pressure}.
Let $(X, T)$ be a dynamical system with a continuous function $f:X\to \mathbb{R}$.
Motivated by statistical mechanics, Ruelle \cite{Ruelle} (in some special cases) and Walters \cite{Walters} (for general systems)
introduced the \textbf{topological pressure} $P(T, f)$ and proved the variational principle 
\begin{equation} \label{eq: standard variational principle for pressure}
  P(T, f) = \sup_{\mu\in \mathcal{M}^T(X)} \left(h_\mu(T) + \int_X f\, d\mu\right). 
\end{equation}

The above (\ref{eq: standard variational principle}) and (\ref{eq: standard variational principle for pressure}) 
are classical and standard in ergodic theory.
Recently, Feng--Huang \cite{Feng--Huang} found an ingenious generalization of this classical theory.
Motivated by fractal geometry of self-affine carpets and sponges \cite{Bedford, McMullen, Kenyon--Peres},
they introduced \textit{weighted versions} of entropy and pressure.

Let $(X, T)$ and $(Y, S)$ be dynamical systems.
A map $\pi:X\to Y$ is called a \textbf{factor map} if $\pi$ is a continuous surjection with 
$\pi \circ T = S\circ \pi$.
We sometimes write $\pi:(X, T)\to (Y, S)$ for clarifying the maps $T$ and $S$.
For an invariant probability measure $\mu\in \mathscr{M}^T(X)$, we denote by $\pi_*\mu \in \mathscr{M}^S(Y)$ 
the push-forward\footnote{This is defined by $\pi_*\mu(A) = \mu\left(\pi^{-1}A\right)$ for $A\subset Y$.} of $\mu$ by $\pi$.
Let $f:X\to \mathbb{R}$ be a continuous function, and let $a_1, a_2$ be two real numbers with $a_1>0$ and $a_2\geq 0$.
Feng--Huang \cite[Question 1.1]{Feng--Huang} asked (and then solved) the following question.

\begin{question} \label{question: main question}
How can one define a meaningful term $P^{(a_1, a_2)}(T, f)$ such that the following variational principle holds?
\[  P^{(a_1, a_2)}(T, f) = \sup_{\mu\in \mathscr{M}^T(X)} \left(a_1 h_\mu(T) + a_2 h_{\pi_*\mu}(S) + \int_X f\, d\mu\right). \]
\end{question}

We describe their approach below. 
It is a modification of the definition of topological entropy given by Bowen \cite{Bowen}, which is in turn 
a modification of the standard definition of Hausdorff dimension.

Here we explain only the case of $f\equiv 0$ for simplicity of the exposition.
For the case of $f\not \equiv 0$, see their paper \cite[\S 3.1]{Feng--Huang}\footnote{They also studied the case that 
a sequence of factor maps $\pi_i:X_i\to X_{i+1}$ ($i=1, 2, \dots, k$) is given. We think that our new approach can be 
also generalized to this setting. But we concentrate on the simplest case in this paper.}.

Let $d$ and $d^\prime$ be metrics on $X$ and $Y$ respectively.
For $x\in X$, a natural number $n$ and $\varepsilon>0$, we define $B_n^{(a_1, a_2)}(x,\varepsilon) \subset X$ as the 
set of $y\in X$ satisfying the following two conditions:
\[  d(T^j x, T^j y) < \varepsilon, \quad (0\leq j < \lceil a_1 n \rceil), \]
\[  d^\prime\left(S^j \pi(x), S^j \pi(x)\right) < \varepsilon, \quad (0\leq j < \lceil (a_1+a_2)n \rceil). \]
Here $\lceil u\rceil$ denotes the least integer not less than $u$.
We call $B_n^{(a_1, a_2)}(x,\varepsilon)$ an \textbf{$(a_1,a_2)$-weighted Bowen ball}.

Let $N$ be a natural number.
We consider families of $(a_1,a_2)$-weighted Bowen balls $\left\{B^{(a_1, a_2)}_{n_j}(x_j, \varepsilon)\right\}_{j=1}^\infty$
satisfying 
\begin{equation}  \label{eq: cover of weighted Bowen balls}
  X = \bigcup_j B^{(a_1,a_2)}_{n_j}(x_j, \varepsilon), \quad n_j \geq N \> (\forall j\geq 1). 
\end{equation}
Let $s\geq 0$. We define $\Lambda^{(a_1, a_2), s}_{N,\varepsilon}(X)$ as the infimum of 
\[ \sum_j \exp(-sn_j) \]
where the infimum is taken over all families $\left\{B^{(a_1, a_2)}_{n_j}(x_j, \varepsilon)\right\}_{j=1}^\infty$
satisfying the above (\ref{eq: cover of weighted Bowen balls}).

The quantity $\Lambda^{(a_1, a_2), s}_{N,\varepsilon}(X)$ is monotone in $N$. So we define 
\[  \Lambda^{(a_1, a_2), s}_\varepsilon(X) = \lim_{N\to \infty} \Lambda^{(a_1, a_2), s}_{N,\varepsilon}(X). \]
We vary the parameter $s$ from $0$ to $\infty$.
There exists a unique value of $s$, which we denote by $h_{\mathrm{top}}^{(a_1, a_2)}(T, \varepsilon)$, where 
the value of $\Lambda^{(a_1, a_2), s}_\varepsilon(X)$ jumps from $\infty$ to $0$:
\[ \Lambda^{(a_1, a_2),s}_\varepsilon(X) 
    =   \begin{cases} 0, \quad &(s>h_{\mathrm{top}}^{(a_1, a_2)}(T, \varepsilon)) \\
                             \infty, \quad &(s< h_{\mathrm{top}}^{(a_1, a_2)}(T, \varepsilon))
         \end{cases}. \]                                                                      
$h_{\mathrm{top}}^{(a_1, a_2)}(T, \varepsilon)$ is monotone in $\varepsilon$. So we define the 
\textbf{$(a_1, a_2)$-weighted topological entropy} of $\pi: X\to Y$ by 
\[ h^{(a_1, a_2)}_{\mathrm{top}}(\pi, T) =  \lim_{\varepsilon\to 0}  h_{\mathrm{top}}^{(a_1, a_2)}(T, \varepsilon). \]

Feng--Huang \cite[Theorem 1.4, Corollary 1.5]{Feng--Huang} solved
Question \ref{question: main question} by this quantity:

\begin{theorem}[Feng--Huang, 2016] \label{theorem: Feng--Huang}
\[  h^{(a_1, a_2)}_{\mathrm{top}}(\pi, T) = \sup_{\mu\in \mathscr{M}^T(X)} \left(a_1 h_\mu(T) + a_2 h_{\pi_*\mu}(S)\right). \]
\end{theorem}

\subsection{New approach}  \label{subsection: new approach}

In the previous subsection we describe the definition of weighted topological entropy 
introduced by Feng--Huang \cite{Feng--Huang}.
In this subsection we describe our new approach.
Our approach is a modification of the familiar definition of topological entropy
(not the Hausdorff-dimension like definition of \cite{Bowen}).

First of all, notice that we can assume $a_1+a_2 = 1$ in Question \ref{question: main question}
because we can reduce the general case to this special case by a simple rescaling.
So we study only this case.
As in the previous subsection, here we explain the entropy case (i.e. the case of $f\equiv 0$) for simplicity.
We will explain the pressure case in \S \ref{section: weighted topological pressure}.

Let $(X,T)$ and $(Y, S)$ be dynamical systems, and let $\pi:X\to Y$ be a factor map.
Let $d$ and $d^\prime$ be metrics on $X$ and $Y$ respectively.
For a natural number $N$ we define metrics $d_N$ and $d^\prime_N$ on $X$ and $Y$ respectively by 
\begin{equation} \label{eq: new metrics}
   d_N(x_1, x_2) = \max_{0\leq n<N} d\left(T^n x_1, T^n x_2\right), \quad 
    d^\prime_N(y_1,y_2) = \max_{0\leq n<N} d^\prime\left(S^n y_1, S^n y_2\right). 
\end{equation}    
For $\varepsilon>0$ and a non-empty subset $\Omega \subset X$ we define 
\begin{equation}  \label{eq: covering number}
 \#\left(\Omega, N, \varepsilon\right) = \min\left\{n\geq 1\middle|
     \parbox{3in}{\centering $\exists$ open subsets $U_1, \dots, U_n$ of $X$ with $\Omega \subset U_1\cup \dots \cup U_n$ and 
     $\diam(U_k, d_N) < \varepsilon$ for all $1\leq k\leq n$}
    \right\}. 
\end{equation}    
Here
$\diam(U_k, d_N) = \sup_{x_1, x_2\in U_k} d_N(x_1, x_2)$ is the diameter of $U_k$ with respect to the metric $d_N$.
When $\Omega$ is the empty set, we define $\#\left(\Omega, N, \varepsilon\right)=0$.
As is well-known, the topological entropy of $(X, T)$ is defined by 
\[ h_{\mathrm{top}}(X,T) = \lim_{\varepsilon \to 0} \left(\lim_{N\to \infty} \frac{\log \#(X, N, \varepsilon)}{N}\right). \]
We will modify this definition.

Let $0\leq w \leq 1$ be a real number. We set 
\begin{equation} \label{eq: weighted covering number}
   \#^w\left(\pi, N, \varepsilon\right) = \min\left\{\sum_{k=1}^n \left(\#\left(\pi^{-1}(V_k), N, \varepsilon\right)\right)^w\middle|
    \parbox{3in}{\centering $Y=V_1\cup \dots \cup V_n$ is an open cover with $\diam \left(V_k, d^\prime_N\right) < \varepsilon$
      for all $1\leq k\leq n$}\right\}. 
\end{equation}      
It is easy to check that this quantity is sub-multiplicative in $N$ and monotone in $\varepsilon$.
So we define the \textbf{$w$-weighted topological entropy} of $\pi:X\to Y$ by
\[ h_{\mathrm{top}}^w(\pi, T)   = \lim_{\varepsilon \to 0} \left(\lim_{N\to \infty} 
    \frac{\log \#^w\left(\pi, N, \varepsilon\right)}{N}\right). \]
This definition uses the metrics $d$ and $d^\prime$, 
but the value of $h_{\mathrm{top}}^w(\pi, T)$ is a topological invariant (i.e. independent of 
the choice of metrics).

The quantity $h^w_{\mathrm{top}}(\pi, T)$ provides another solution to Question \ref{question: main question}
for the case of $f\equiv 0$ and $(a_1, a_2) = (w, 1-w)$.
This is our main result for the weighted topological entropy.

\begin{theorem}[\textbf{Variational principle for $w$-weighted topological entropy}] \label{theorem: main theorem on p-weighted entropy}
For $0\leq w\leq 1$
\[ \h^w(\pi, T) = \sup_{\mu\in \mathscr{M}^T(X)} \left\{w h_\mu(T) + (1-w)h_{\pi_*\mu}(S)\right\}. \]
\end{theorem}

As the above definition of $\h^w(\pi, T)$ is close to the standard definition of topological entropy,
the proof of this theorem is also close to a well-known proof of the standard variational principle.
The basic structure of the proof is borrowed from the famous argument
of Misiurewicz \cite{Misiurewicz}.
At some technical points, we use the theory of \textit{principal extensions} 
\cite{Downarowicz, Downarowicz--Huczek}.

By combining Theorems \ref{theorem: Feng--Huang} and \ref{theorem: main theorem on p-weighted entropy}, we get a corollary:

\begin{corollary}  \label{corollary: Feng--Huang = ours}
\[  \h^{(w, 1-w)}(\pi, T) = \h^w(\pi, T). \]
Here the left-hand side is the weighted topological entropy $\h^{(a_1, a_2)}(\pi, T)$
for $(a_1, a_2) = (w, 1-w)$ defined in the previous subsection.
\end{corollary} 

This corollary seems to be a very interesting statement.
The author cannot see any direct way to prove it (without using the variational principles).

\begin{problem}
Can one prove the equality $\h^{(w, 1-w)}(\pi, T) = \h^w(\pi, T)$ without using measure theory?
\end{problem}

The following example illustrates the importance of the equality $\h^{(w, 1-w)}(\pi, T) = \h^w(\pi, T)$.

\begin{example}[Bedford--McMullen carpets] \label{example: Bedford--McMullen carpets}
Let $\mathbb{T} = \mathbb{R}/\mathbb{Z}$ be the circle, and let $\mathbb{T}^2 = \mathbb{R}^2/\mathbb{Z}^2$ be the torus.
Let $a$ and $b$ be two natural numbers with $a\geq b\geq 2$.  
Set $A = \{0,1,2,\dots, a-1\}$ and $B = \{0,1,2,\dots, b-1\}$.
Let $R\subset A\times B$ be a non-empty subset, and define
\[  R^\prime = \{y\in B|\, (x, y) \in R \text{ for some $x\in A$}\}. \]
We define $X\subset \mathbb{T}^2$ and $Y\subset \mathbb{T}$ by
\[ X := \left\{\left(\sum_{n=1}^\infty \frac{x_n}{a^n}, \sum_{n=1}^\infty \frac{y_n}{b^n}\right) \in \mathbb{T}^2 \middle|\, 
             (x_n, y_n)\in R \text{ for all $n\geq 1$}\right\}, \]
\[ Y := \left\{\sum_{n=1}^\infty \frac{y_n}{b^n} \in \mathbb{T}\middle|\,   y_n \in R^\prime \text{ for all $n\geq 1$}\right\}. \]
The space $X$ is the famous Bedford--McMullen carpet \cite{Bedford, McMullen}.
We are going to explain that we can calculate the Hausdorff dimension of $X$ (with respect to the natural metric on $\mathbb{T}^2$)
by using Corollary \ref{corollary: Feng--Huang = ours}.

We define continuous maps $T:X\to X$ and $S:Y\to Y$ by 
\[ T(x, y) = (ax, by), \quad S(y) = by. \]
$(X, T)$ and $(Y, S)$ are dynamical systems.
Let $\pi:X\to Y$ be the natural projection.
$\pi$ is a factor map between $(X, T)$ and $(Y, S)$.
We are interested in its weighted topological entropy.
Set 
\[ w = \frac{\log b}{\log a} = \log_a b. \] 
We have $0\leq w \leq 1$.
It directly follows from the definitions\footnote{The $(a_1, a_2)$-weighted Bowen ball $B_n^{(a_1, a_2)}(x,\varepsilon)$
for $a_1 = \log_a b$ and $a_2=1-\log_a b$ is approximately a square of side length $b^n$.}
in \S \ref{subsection: weighted topological entropy and pressure} that 
the Hausdorff dimension of $X$ is given by 
\[ \dim_H X = \frac{\h^{(w, 1-w)}(\pi, T)}{\log b}. \]
From the equality $\h^{(w, 1-w)}(\pi, T) = \h^w(\pi, T)$ in Corollary \ref{corollary: Feng--Huang = ours}, 
we also have 
\begin{equation} \label{eq: hausdorff dimension and weighted topological entropy}
   \dim_H X = \frac{\h^w(\pi, T)}{\log b}. 
\end{equation}   

Now we calculate the $w$-weighted topological entropy $h^{w}(\pi, T)$:

\begin{claim} \label{claim: calculation of h^w}
For each $y\in B$ we define $t(y)$ as the number of $x\in A$ satisfying $(x, y)\in R$.
Then 
\[  \h^{w}(\pi, T) = \log\left(\sum_{y\in R^\prime} t(y)^w\right). \]
\end{claim}

\begin{proof}
First notice that, in the definitions (\ref{eq: covering number}) and (\ref{eq: weighted covering number}), we can use 
closed covers instead of open covers; this does not change their values.
Here we will consider closed covers.

We define a metric $d^\prime$ on $\mathbb{T}$ by 
\[ d^\prime(x_1, x_2) = \min_{n\in \mathbb{Z}} |x_1-x_2-n|. \]
We define a metric $d$ on $\mathbb{T}^2$ by 
\[ d\left((x_1, y_1), (x_2, y_2)\right) = \max\left(d^\prime(x_1, x_2), d^\prime(y_1, y_2)\right). \]

Let $\varepsilon>0$ and take a natural number $m$ with $b^{-m} < \varepsilon$.
Let $N$ be a natural number. For each $v\in (R^\prime)^{N+m}$, set 
\[ V_v = \left\{\sum_{n=1}^\infty \frac{y_n}{b^n} \in Y\middle|\, (y_1, \dots, y_{N+m}) = v \right\}. \]
These form a closed covering of $Y$ with $\diam\left(V_v, d^\prime_N\right) < \varepsilon$.
For each $(u, v)\in R^{N+m} \subset A^{N+m}\times B^{N+m}$ (where $u\in A^{N+m}$ and $v\in (R^\prime)^{N+m}$),
we set 
\[ U_{(u, v)} = \left\{\left(\sum_{n=1}^\infty \frac{x_n}{a^n}, \sum_{n=1}^\infty \frac{y_n}{b^n}\right) \in X\middle|\, 
                    (x_1, \dots, x_{N+m}) = u, \, (y_1, \dots, y_{N+m}) = v\right\}. \]
These are closed subsets of $X$ with $\diam\left(U_{(u,v)}, d_N\right) < \varepsilon$ and 
\[ \pi^{-1}(V_v) = \bigcup_{\substack{u\in A^{N+m} \\ \text{with $(u,v)\in R^{N+m}$}}} U_{(u, v)}. \]
Hence, for $v=(v_1, \dots, v_{N+m}) \in (R^\prime)^{N+m}$ 
\[  \#\left(\pi^{-1}(V_v), N, \varepsilon\right) \leq t(v_1)\cdots t(v_{N+m}). \]
Therefore 
\[ \#^w\left(\pi, N, \varepsilon\right) \leq \sum_{v_1, \dots, v_{N+m}\in R^\prime} \left(t(v_1)\cdots t(v_{N+m})\right)^{w}
     = \left(\sum_{v\in R^\prime} t(v)^w\right)^{N+m}. \]
Thus 
\[  h_{\mathrm{top}}^w(\pi, T)   = \lim_{\varepsilon \to 0} \left(\lim_{N\to \infty} 
    \frac{\log \#^w\left(\pi, N, \varepsilon\right)}{N}\right)
     \leq \log\left(\sum_{y\in B} t(y)^w\right). \]

Next, let $0<\varepsilon <\frac{1}{a}$.
Fix $(p, q) \in R$.
For a natural number $N$, we consider the following points in $Y$:
\begin{equation}  \label{eq: separated set in Y}
  \sum_{n=1}^N \frac{v_n}{b^n} + \sum_{n=N+1}^\infty \frac{q}{b^n}, \quad (v_1, \dots, v_N\in R^\prime). 
\end{equation}  
These points form an $\varepsilon$-separated set in $Y$ with respect to the metric $d^\prime_N$. 
We also consider the following points in $X$:
\begin{equation}  \label{eq: separated set in X}
 \left(\sum_{n=1}^N \frac{u_n}{a^n}, \sum_{n=1}^N \frac{v_n}{b^n} \right)
        + \sum_{n=N+1}^\infty \left(\frac{p}{a^n},  \frac{q}{b^n}\right), \quad 
    \left((u_1, v_1), \dots, (u_N, v_N)\in R\right). 
\end{equation}    
These points form an $\varepsilon$-separated set in $X$ with respect to the metric $d_N$.

Suppose $Y=V_1\cup \dots \cup V_n$ is a covering with $\diam(V_k, d^\prime_N) < \varepsilon$.
Then each $V_k$ contains at most one point of (\ref{eq: separated set in Y}).
If $V_k$ contains a point $\sum_{n=1}^N \frac{v_n}{b^n} + \sum_{n=N+1}^\infty \frac{q}{b^n}$, 
then $\pi^{-1}(V_k)$ contains 
$t(v_1)\cdots t(v_N)$ points of the form (\ref{eq: separated set in X}) and hence 
\[ \#\left(\pi^{-1}(V_k), N, \varepsilon\right) \geq t(v_1) \cdots t(v_N). \]
So 
\[ \#^w\left(\pi, N, \varepsilon\right) \geq \sum_{v_1, \dots, v_{N}\in R^\prime} \left(t(v_1)\cdots t(v_{N})\right)^{w}
     = \left(\sum_{v\in R^\prime} t(v)^w\right)^{N}. \]
This shows
\[  h_{\mathrm{top}}^w(\pi, T)   = \lim_{\varepsilon \to 0} \left(\lim_{N\to \infty} 
    \frac{\log \#^w\left(\pi, N, \varepsilon\right)}{N}\right)
     \geq \log\left(\sum_{y\in R^\prime} t(y)^w\right). \]
Notice that this proof of the claim is completely elementary.
We have not used any sophisticated technique (in particular, measure theory).
\end{proof}

From (\ref{eq: hausdorff dimension and weighted topological entropy}) and Claim \ref{claim: calculation of h^w}
\begin{equation}  \label{eq: dimension of Bedford--McMullen carpet}
 \dim_H X = \frac{\log\left(\sum_{y\in R^\prime} t(y)^w\right)}{\log b}  
    = \log_b \left(\sum_{y\in R^\prime} t(y)^{\log_a b}\right). 
\end{equation}   
This is a famous formula for the Hausdorff dimension of the Bedford--McMullen carpet \cite{Bedford, McMullen}.
Therefore we conclude that the equality $\h^{(1-w,w)}(\pi, T) = \h^w(\pi, T)$ provides this famous formula
fairly easily.
This suggests that the equality $\h^{(1-w, w)}(\pi, T) = \h^w(\pi, T)$ is a rather deep statement.
We can say that it is a \textit{topological generalization} of the dimension formula for the Bedford--McMullen carpet.

Kenyon--Peres \cite[Theorem 1.1, Theorem 3.2]{Kenyon--Peres sofic} generalized 
the formula (\ref{eq: dimension of Bedford--McMullen carpet}) to 
closed $T$-invariant subsets of $\mathbb{T}^2$ which correspond to subshifts of finite type or sofic subshifts under the natural Markov partition.
We can also prove their results from the equality $\h^{(1-w,w)}(\pi, T) = \h^w(\pi, T)$ as in the above.
\end{example}

The above example also illustrates that the two notions $\h^{(a_1, a_2)}(\pi, T)$ and 
$\h^w(\pi, T)$ have their own advantages.
One of the great advantages of $\h^{(a_1, a_2)}(\pi, T)$ is that its definition is intrinsically
related to Hausdorff dimension.
So it can be directly applied to the study of geometric measure theory.
On the other hand, the advantage of $\h^w(\pi, T)$ is that its definition is elementary and hence (sometimes) easy to calculate.

In \cite[p. 441]{Feng--Huang}, Feng--Huang asked how to generalize their result to $\mathbb{Z}^d$-actions.
It seems rather straightforward to generalize our new approach to $\mathbb{Z}^d$-actions 
and, possibly, actions of amenable groups.

\begin{problem}
Suppose that both $\h^{(a_1, a_2)}(\pi, T)$ and $\h^w(\pi, T)$ are generalized to group actions.
Can one deduce any interesting consequence of their coincidence 
(like the above calculation of the Hausdorff dimension of the Bedford--McMullen carpet)?
\end{problem}

We would like to mention the papers of Barral--Feng \cite{Barral--Feng arXiv, Barral--Feng} and Feng \cite{Feng} 
(see also Yayama \cite{Yayama compensation, Yayama}).
These papers studied Question \ref{question: main question} and related questions when $(X, T)$ and $(Y, S)$ are 
subshifts over finite alphabets.
When $(X, T)$ and $(Y, S)$ are subshifts, the above definition of $\h^w(\pi, T)$ 
(and its pressure version in \S \ref{section: weighted topological pressure}) is essentially the same with 
the one given in \cite[Theorem 1.1]{Barral--Feng arXiv} (see also \cite[Theorem 3.1]{Barral--Feng}).
So we can say that the above definition generalizes the approach in \cite[Theorem 1.1]{Barral--Feng arXiv} 
from subshifts to general dynamical systems.

This paper studies only the abstract theory of $\h^w(\pi, T)$ and its pressure version.
But the main motivation for the author to introduce these quantities is not to develop the abstract theory.
The author naturally came up with the above definition of $\h^w(\pi, T)$ when he studied mean Hausdorff dimension of
certain \textit{infinite dimensional fractals}.
(Mean Hausdorff dimension is a dynamical version of Hausdorff dimension introduced in \cite{Lindenstrauss--Tsukamoto double VP}.)
We plan to describe this connection in a separate paper.

\section{Weighted topological pressure}  \label{section: weighted topological pressure}

In this section we introduce our new definition of weighted topological pressure.
For the original approach, see \cite[\S 3.1]{Feng--Huang}.

Let $\pi:X\to Y$ be a factor map from a dynamical systems $(X, T)$ to a dynamical system $(Y, S)$.
Let $f:X\to \mathbb{R}$ be a continuous function.

Let $d$ and $d^\prime$ be metrics on $X$ and $Y$ respectively.
For a natural number $N$ we define new metrics $d_N$ and $d^\prime_N$ on $X$ and $Y$ respectively 
by (\ref{eq: new metrics}). 
We also define a continuous function $\mathbb{S}_N f:X\to \mathbb{R}$ by 
\[ \mathbb{S}_N f(x) = f(x) + f(Tx) + f(T^2 x) + \dots + f(T^{N-1}x). \]
The metrics $d_N, d^\prime_N$ and function $\mathbb{S}_N f$ are sometimes denoted by 
$d^T_N, (d^\prime)^S_N$ and $\mathbb{S}^T_N f$ respectively for clarifying the underlying dynamics.

For $\varepsilon>0$ and a non-empty subset $\Omega\subset X$ we define 
\begin{equation} \label{eq: pressure}
  P\left(\Omega, f, N, \varepsilon\right) = \inf\left\{\sum_{k=1}^n \exp\left(\sup_{U_k} \mathbb{S}_N f\right) \middle|
     \parbox{3in}{\centering $\exists$ open subsets $U_1, \dots, U_n$ of $X$ with $\Omega \subset U_1\cup \dots \cup U_n$ and 
     $\diam(U_k, d_N) < \varepsilon$ for all $1\leq k\leq n$}
    \right\}. 
\end{equation}    
(When $U_k$ is the empty set, we assume that the term $\exp \left(\sup_{U_k} \mathbb{S}_N f\right)$ is zero.)
We sometimes denote $P\left(\Omega, f, N, \varepsilon\right)$ by $P_T\left(\Omega, f, N, \varepsilon\right)$
for clarifying the map $T$.
When $\Omega$ is the empty set, we define $P\left(\Omega, f, N, \varepsilon\right) = 0$.
It is well-known that the topological pressure of $(X, T, f)$ is given by 
\[ P(T, f) = \lim_{\varepsilon\to 0} \left(\lim_{N\to \infty} \frac{\log P\left(X, f, N, \varepsilon\right)}{N}\right). \]
We will modify this definition.
Let $0\leq w \leq 1$ be a real number.
We set 
\begin{equation}  \label{eq: weighted pressure in step N}
   \begin{split}
    & P^w\left(\pi, f, N,\varepsilon\right)  \\
    & =  \inf\left\{\sum_{k=1}^n \left(P\left(\pi^{-1}(V_k), f, N, \varepsilon\right)\right)^w\middle|\, 
       \parbox{3in}{\centering $Y=V_1\cup \dots \cup V_n$ is an open cover with $\diam \left(V_k, d^\prime_N\right) < \varepsilon$
      for all $1\leq k\leq n$}\right\}. 
    \end{split}
\end{equation}
We sometimes denote this by $P^w_T\left(\pi, f, N, \varepsilon\right)$.

The quantity $P^w\left(\pi, f, N,\varepsilon\right)$ is sub-multiplicative in $N$ and monotone in $\varepsilon$.
So we define the \textbf{$w$-weighted topological pressure} by 
\[  P^w\left(\pi, T, f\right) = \lim_{\varepsilon\to 0} \left(\lim_{N\to \infty} \frac{\log P^w\left(\pi, f, N,\varepsilon\right)}{N}\right). \]
The value of $P^w\left(\pi, T, f\right)$ is independent of the choices of the metrics $d$ and $d^\prime$.
So it provides a topological invariant.
We sometimes use the notation $P^w\left(\pi, X, T, Y, S, f\right)$ instead of $P^w\left(\pi, T, f\right)$
for clarifying all the data involved.

Now we state our main result of the paper.

\begin{theorem}[\textbf{Variational principle for $w$-weighted topological pressure}]  \label{theorem: main theorem}
For any $0\leq w\leq 1$
\[    P^w\left(\pi, T, f\right) = \sup_{\mu\in \mathscr{M}^T(X)} \left(w h_\mu(T) + (1-w)h_{\pi_*\mu}(S) + w\int_X f\, d\mu\right). \]
\end{theorem}

When $f\equiv 0$, we have $P^w\left(\pi, T, f\right) = \h^w\left(\pi, T\right)$. 
So Theorem \ref{theorem: main theorem on p-weighted entropy} in \S \ref{subsection: new approach}
follows from Theorem \ref{theorem: main theorem}.
The proof of Theorem \ref{theorem: main theorem} occupies all the rest of the paper.

For the simplicity of the notation, we write 
\begin{equation}  \label{eq: weighted pressure via measures}
    P^w_{\mathrm{var}}\left(\pi, T, f\right) := 
   \sup_{\mu\in \mathscr{M}^T(X)} \left(w h_\mu(T) + (1-w)h_{\pi_*\mu}(S) + w\int_X f\, d\mu\right). 
\end{equation}   
(Here \textit{var} is the abbreviation of \textit{variational}.) Then our main purpose is to prove the equality 
\[ P^w\left(\pi, T, f\right)   =  P^w_{\mathrm{var}}\left(\pi, T, f\right). \]

In the rest of this section we gather some elementary properties of $w$-weighted topological pressure.
Here we always assume that $\pi:(X, T)\to (Y, S)$ is a factor map between dynamical systems
with a continuous function $f:X\to \mathbb{R}$.
We take $0\leq w\leq 1$.
Let $d$ and $d^\prime$ be metrics on $X$ and $Y$ respectively.

\begin{lemma} \label{lemma: amplification}
Let $m$ be a natural number.
\[  P^w\left(\pi, T^m, \mathbb{S}_m f\right) = m P^w\left(\pi, T, f\right). \]
Here the left-hand side is $P^w\left(\pi, X,T^m,Y,S^m, \mathbb{S}^T_m f\right)$.
\end{lemma}

\begin{proof}
Let $\varepsilon$ be a positive number.
There exists $0<\delta<\varepsilon$ such that 
\[ d(x_1, x_2) < \delta \Longrightarrow d^T_m(x_1, x_2) < \varepsilon, \quad (x_1, x_2\in X), \]
\[ d^\prime(y_1, y_2) < \delta \Longrightarrow (d^\prime)^S_m(y_1,y_2) < \varepsilon, \quad 
    (y_1, y_2\in Y).  \]
Then for any natural number $N$
\[ d^{T^m}_N(x_1, x_2) < \delta \Longrightarrow d^T_{mN}(x_1, x_2) < \varepsilon, \quad 
    (x_1, x_2\in X),  \]
\[ (d^\prime)^{S^m}_N(y_1,y_2) < \delta \Longrightarrow (d^\prime)^S_{mN}(y_1,y_2) < \varepsilon, \quad 
    (y_1, y_2\in Y).  \]
Since $\mathbb{S}^{T^m}_N \left(\mathbb{S}^T_m f\right) = \mathbb{S}^T_{mN} f$, for any subset $\Omega \subset X$
\[  P_{T^m}\left(\Omega, \mathbb{S}^T_m f, N,\varepsilon\right) \leq 
     P_T\left(\Omega, f, mN, \varepsilon\right) \leq 
     P_{T^m}\left(\Omega, \mathbb{S}^T_m f, N, \delta\right). \] 
Then 
\[ P^w_{T^m}\left(\pi, \mathbb{S}^T_m f, N,\varepsilon\right) \leq 
   P^w_T\left(\pi, f, mN, \varepsilon\right)  \leq 
   P^w_{T^m}\left(\pi, \mathbb{S}^T_m f, N, \delta \right). \]
Thus 
\[  P^w\left(\pi, T^m, \mathbb{S}^T_m f\right) = m P^w\left(\pi, T, f\right). \]
\end{proof}

\begin{lemma}  \label{lemma: pressure and factor map}
Let $(X^\prime, T^\prime)$ be a dynamical system, and let $\varphi: (X^\prime, T^\prime) \to (X, T)$ be a factor map.
\[ \xymatrix{
   (X^\prime, T^\prime)  \ar[rr]^{\varphi} \ar[rrd]_{\pi\circ \varphi} & & (X, T)  \ar[d]^{\pi} \\
   &  & (Y,S)}
 \]
Then 
\[  P^w(\pi, T, f) \leq  P^w\left(\pi\circ \varphi, T^\prime, f\circ \varphi\right). \]
Here the right-hand side is $P^w\left(\pi\circ \varphi, X^\prime, T^\prime, Y, S, f\circ \varphi\right)$.
\end{lemma}

\begin{proof}
Let $\tilde{d}$ be a metric on $X^\prime$.
For any $\varepsilon>0$ there exists $0<\delta<\varepsilon$ satisfying 
\[ \tilde{d}(x_1, x_2) < \delta \Longrightarrow d(x_1, x_2) < \varepsilon. \]
Then for any $N>0$
\[ \tilde{d}_N(x_1, x_2) < \delta \Longrightarrow d_N(x_1, x_2) < \varepsilon. \]
From this, we have for any $\Omega \subset X^\prime$
\[  P_T\left(\varphi(\Omega), f, N, \varepsilon\right)  \leq  P_{T^\prime}\left(\Omega, f\circ \varphi, N, \delta\right). \]
For any $V\subset Y$
\[  \varphi\left((\pi\circ \varphi)^{-1}(V)\right) = \pi^{-1}(V). \]
So 
\[ P_T\left(\pi^{-1}(V), f, N, \varepsilon\right) \leq 
    P_{T^\prime}\left((\pi\circ \varphi)^{-1}(V), f \circ \varphi, N, \delta\right). \]

Then 
\[ P^w_T\left(\pi, f, N,\varepsilon\right)  \leq P^w_{T^\prime}\left(\pi\circ \varphi, f\circ \varphi, N, \delta \right). \]
Therefore 
\[  P^w(\pi, T, f) \leq  P^w\left(\pi\circ \varphi, T^\prime, f\circ \varphi\right). \]
\end{proof}

The next lemma is a bit complicated.
It might be better for some readers to look at Remark \ref{remark: weighted pressure and factor} below before reading 
the lemma. It will provide a clearer perspective.

\begin{lemma} \label{lemma: fiber product}
Let $(Y^\prime, S^\prime)$ be a dynamical system, and let $\phi:(Y^\prime, S^\prime)\to (Y,S)$ be a factor map.
Define the \textbf{fiber product}
\[ X\times_Y Y^\prime = \left\{(x, y)\in X\times Y^\prime\middle|\, \pi(x) = \phi(y)\right\}. \]
$\left(X\times_Y Y^\prime, T\times S^\prime\right)$ becomes a dynamical system.
We define factor maps $\varphi: X\times_Y Y^\prime\to X$ and $\Pi:X\times_Y Y^\prime\to Y^\prime$ by 
\[ \varphi(x, y) = x, \quad \Pi(x,y) = y. \]
The diagram is as follows:
\[  \xymatrix{ \left(X\times_Y Y^\prime, T\times S^\prime\right) \ar[rrr]^{\varphi} \ar[d]_{\Pi} & & & (X,T) \ar[d]^{\pi} \\
                       (Y^\prime, S^\prime) \ar[rrr]_{\phi} & & & (Y,S) } \]
Then 
\[  P^w\left(\pi, T, f\right)  \leq  P^w\left(\Pi, T\times S^\prime, f\circ \varphi\right). \]
Here the right-hand side is $P^w\left(\Pi, X\times_Y Y^\prime, T\times S^\prime, Y^\prime, S^\prime, f\circ \varphi\right)$.
\end{lemma}

\begin{proof}
The point of the proof is that for any subset $A\subset Y^\prime$ we have
\[  \pi^{-1}\left(\phi(A)\right) = \varphi\left(\Pi^{-1}(A)\right). \]

Let $\tilde{d}$ be a metric on $Y^\prime$ and we define a metric $\rho$ on $X\times_Y Y^\prime$ by 
\[ \rho\left((x_1, y_1), (x_2, y_2)\right) = \max\left(d(x_1,x_2), \tilde{d}(y_1,y_2)\right). \]
Let $\varepsilon$ be a positive number. We have 
\[ \rho\left((x_1, y_1), (x_2, y_2)\right)  < \varepsilon \Longrightarrow d(x_1, x_2) < \varepsilon. \]
Then for any natural number $N$ and any subset $\Omega \subset X\times_Y Y^\prime$
\[  P_T\left(\varphi(\Omega), f, N, \varepsilon\right) \leq 
     P_{T\times S^\prime}\left(\Omega, f\circ \varphi, N, \varepsilon\right). \]
In particular, for any subset $A \subset Y^\prime$
\begin{equation}  \label{eq: fiber product and weighted pressure key estimate}
  \begin{split}
    P_T\left(\pi^{-1}\left(\phi(A)\right), f, N, \varepsilon\right) 
    &  = P_T\left(\varphi\left(\Pi^{-1}(A)\right), f, N, \varepsilon\right)   \\
     &  \leq P_{T\times S^\prime}\left(\Pi^{-1}(A), f\circ \varphi, N, \varepsilon\right). 
  \end{split}   
\end{equation}

There exists $0<\delta<\varepsilon$ such that 
\[  \tilde{d}(y_1, y_2) < \delta \Longrightarrow d^\prime\left(\phi(y_1), \phi(y_2)\right) < \varepsilon. \]
Now we claim that 
\[ P^w_T\left(\pi, f, N, \varepsilon\right) \leq  P^w_{T\times S^\prime}\left(\Pi, f\circ \varphi, N,\delta\right). \]
Indeed take any positive number $C$ with 
\[ P^w_{T\times S^\prime}\left(\Pi, f\circ \varphi, N,\delta\right)  < C. \]
Then there exists an open covering $Y^\prime = V_1\cup \dots \cup V_n$ such that 
$\diam\left(V_k, \tilde{d}_N\right) < \delta$ for all $1\leq k \leq n$ and 
\[ \sum_{k=1}^n \left(P_{T\times S^\prime}\left(\Pi^{-1}(V_k), f\circ \varphi, N, \delta \right)\right)^w  < C. \]
We can find compact subsets $A_k \subset V_k$ satisfying $Y^\prime = A_1\cup \dots \cup A_n$.
We have 
\begin{equation*}
  \begin{split}
    \sum_{k=1}^n \left(P_T\left(\pi^{-1}\left(\phi(A_k)\right), f, N, \varepsilon\right)\right)^w 
     & \leq \sum_{k=1}^n \left(P_{T\times S^\prime}\left(\Pi^{-1}(A_k), f\circ \varphi, N, \varepsilon \right)\right)^w  
     \quad \text{by (\ref{eq: fiber product and weighted pressure key estimate})} \\
     & \leq \sum_{k=1}^n \left(P_{T\times S^\prime}\left(\Pi^{-1}(A_k), f\circ \varphi, N, \delta \right)\right)^w
      \quad \text{by $\delta<\varepsilon$}  \\
     &\leq \sum_{k=1}^n \left(P_{T\times S^\prime}\left(\Pi^{-1}(V_k), f\circ \varphi, N, \delta \right)\right)^w 
     \quad \text{by $A_k\subset V_k$} \\
     &      < C.
  \end{split}
\end{equation*}     
Each $\phi(A_k)$ is a closed subset of $Y$ with
$\diam\left(\phi(A_k), d^\prime_N\right) < \varepsilon$.
By the definition (\ref{eq: pressure}),
there exist open subsets $W_k \supset \phi(A_k)$ of $Y$ for $1\leq k\leq n$ such that 
$\diam\left(W_k, d^\prime_N\right) < \varepsilon$ and 
\[  \sum_{k=1}^n \left(P\left(\pi^{-1}\left(W_k\right), f, N, \varepsilon\right)\right)^w  < C. \]
Noticing $Y = W_1\cup \dots \cup W_n$, we have
\[ P^w_T\left(\pi, f, N, \varepsilon\right)  <  C. \]
Since $C$ is an arbitrary number larger than $P^w_{T\times S^\prime}\left(\Pi, f\circ \varphi, N,\delta\right)$, 
this shows
\[  P^w_T\left(\pi, f, N, \varepsilon\right)  \leq   P^w_{T\times S^\prime}\left(\Pi, f\circ \varphi, N,\delta\right). \]
Thus we conclude
\[ P^w\left(\pi, T, f\right)  \leq  P^w\left(\Pi, T\times S^\prime, f\circ \varphi\right). \]
\end{proof}

\begin{remark}  \label{remark: weighted pressure and factor}
Let $(X^\prime, T^\prime)$ and $(Y^\prime, S^\prime)$ be dynamical systems, and let 
$\pi^\prime:X^\prime \to Y^\prime$ be a factor map.
Suppose there exist factor maps $\varphi:(X^\prime, T^\prime) \to (X, T)$ and 
$\phi:(Y^\prime, S^\prime)\to (Y, S)$ satisfying $\pi\circ \varphi = \phi\circ \pi^\prime$.
\[  \xymatrix{ (X^\prime, T^\prime) \ar[r]^{\varphi} \ar[d]_{\pi^\prime} & (X,T) \ar[d]^{\pi} \\
                    (Y^\prime, S^\prime) \ar[r]_{\phi} &  (Y, S) }\]
Then 
\begin{equation}  \label{eq: pressure and factor maps in general}
  P^w\left(\pi, T, f\right)  \leq P^w\left(\pi^\prime, T^\prime, f\circ \varphi\right) 
\end{equation}  
Here the right-hand side is $P^w\left(\pi^\prime, X^\prime, T^\prime, Y^\prime, S^\prime, f\circ \varphi\right)$.
The above Lemmas \ref{lemma: pressure and factor map} and \ref{lemma: fiber product}
are special cases of this statement.
We can prove (\ref{eq: pressure and factor maps in general}) by using the variational principle (Theorem \ref{theorem: main theorem}).
But it seems difficult to prove it in an elementary way.
We will not use (\ref{eq: pressure and factor maps in general}) in the paper.
\end{remark}

Finally we mention two basic results on calculus, which underpin many arguments of this paper:

\begin{lemma}  \label{lemma: calculus}
   \begin{enumerate}
    \item[(1)] For $0\leq w \leq 1$ and nonnegative numbers $x, y$
    \[  (x+y)^w  \leq x^w + y^w. \]
    \item[(2)] Let $p_1, \dots, p_n$ be nonnegative numbers with $p_1+\dots+p_n = 1$.
    For any real numbers $x_1, \dots, x_n$
    \[  \sum_{i=1}^n \left(-p_i \log p_i + p_i x_i\right)  \leq \log \sum_{i=1}^n e^{x_i}. \] 
    In particular (letting $x_1 = \dots = x_n = 0$)
    \[ - \sum_{i=1}^n p_i \log p_i \leq \log n. \]
   \end{enumerate}
\end{lemma}

\begin{proof}
(1) is completely elementary. (2) is proved in \cite[\S 9.3, Lemma 9.9]{Walters book}.
\end{proof}

\section{Kolmogorov--Sinai entropy}  \label{section: Kolmogorov--Sinai entropy}

In this section we review basic definitions on Kolmogorov--Sinai entropy.
For the details, see the book of Walters \cite{Walters book}.

Let $(X, \mu)$ be a probability measure space, namely $X$ is a set equipped with a $\sigma$-algebra and
$\mu$ is a probability measure defined on it.
In our later applications, $X$ is always a compact metrizable space with the standard Borel $\sigma$-algebra.

Let $\mathscr{A} = \{A_1, A_2, \dots, A_n\}$ be a finite measurable partition of $X$, 
namely each $A_i$ is a measurable subset of $X$ and 
\[ X = \bigcup_{i=1}^n A_i, \quad A_i \cap A_j = \emptyset \quad (i\neq j). \]
We define the \textbf{Shannon entropy} of $\mathscr{A}$ by 
\[ H_\mu\left(\mathscr{A}\right) = -\sum_{i=1}^n \mu(A_i) \log \mu(A_i),  \]
where we assume $0 \log 0 = 0$.

For another finite measurable partition $\mathscr{A}^\prime = \{A^\prime_1, A^\prime_2, \dots, A^\prime_m\}$ we set 
\[ \mathscr{A}  \vee \mathscr{A}^\prime 
    = \left\{A_i \cap A^\prime_j \middle|\, 1\leq i \leq n, 1\leq j\leq m\right\}. \]
This is a finite measurable partition of $X$.
We define the \textbf{conditional entropy} by 
\[  H_\mu\left(\mathscr{A}\middle| \mathscr{A}^\prime\right) = -\sum_{\substack{1\leq j\leq m \\ \text{with $\mu(A^\prime_j) > 0$}}}
     \mu(A_j^\prime) \left\{\sum_{i=1}^n 
     \frac{\mu\left(A_i \cap A^\prime_j\right)}{\mu(A^\prime_j)} \log \frac{\mu\left(A_i \cap A^\prime_j\right)}{\mu(A^\prime_j)}\right\}.  \]
Here, in the first summation, we have considered only the index $j$ satisfying $\mu(A^\prime_j) > 0$.
We have (\cite[Theorem 4.3 (i)]{Walters book})
\[ H_\mu\left(\mathscr{A}\vee \mathscr{A}^\prime\right) 
  = H_\mu(\mathscr{A}^\prime) + H_\mu\left(\mathscr{A}\middle| \mathscr{A}^\prime\right). \]

We write $\mathscr{A}^\prime \prec \mathscr{A}$ if $\mathscr{A}\vee \mathscr{A}^\prime = \mathscr{A}$.
This is equivalent to the condition that for every $A \in \mathscr{A}$ there exists $A^\prime \in \mathscr{A}^\prime$
containing $A$. If $\mathscr{A}^\prime \prec \mathscr{A}$ then 
\[  H_\mu\left(\mathscr{A}|\mathscr{A}^\prime\right) = H_\mu(\mathscr{A}) - H_\mu(\mathscr{A}^\prime) \]
and $H_\mu(\mathscr{A}^\prime) \leq H_\mu(\mathscr{A})$.

\begin{lemma} \label{lemma: subadditivity and concavity of entropy}
 \begin{enumerate}
   \item[(1)] $H_\mu(\mathscr{A})$ is subadditive in $\mathscr{A}$. Namely for two finite measurable partitions 
   $\mathscr{A}$ and $\mathscr{A}^\prime$ of $X$
   \[  H_\mu\left(\mathscr{A}\vee \mathscr{A}^\prime\right)  \leq H_\mu\left(\mathscr{A}\right)  + H_\mu\left(\mathscr{A}^\prime\right). \]
   \item[(2)]  $H_\mu(\mathscr{A})$ is concave in $\mu$. 
   Namely for $0\leq t\leq 1$ and two probability measures $\mu$ and $\mu^\prime$ on $X$ 
   \[  H_{(1-t)\mu + t\mu^\prime}(\mathscr{A}) \geq (1-t)H_\mu(\mathscr{A}) + t H_\mu(\mathscr{A}). \]
 \end{enumerate}
\end{lemma}

\begin{proof}
See \cite[Theorem 4.3 (viii)]{Walters book} and \cite[\S 8.1 Remark]{Walters book} for the proofs of 
(1) and (2) respectively.
\end{proof}

Let $T:X\to X$ be a measurable map satisfying $T_*\mu = \mu$.
Let $\mathscr{A}$ be a finite measurable partition of $X$.
For a natural number $N$ 
we define a new measurable partition $\mathscr{A}^N$ of $X$ by 
\[ \mathscr{A}^N = \mathscr{A}\vee T^{-1}\mathscr{A} \vee T^{-2}\mathscr{A} \vee \cdots \vee T^{-(N-1)}\mathscr{A}.  \]
We define the entropy $h_\mu(T, \mathscr{A})$ by 
\[  h_\mu(T, \mathscr{A}) = \lim_{N\to \infty} \frac{H_\mu\left(\mathscr{A}^N\right)}{N}. \]
Finally we define the \textbf{Kolmogorov--Sinai entropy} of the measure-preserving transformation $T$ by 
\[ h_\mu(T) = \sup\left\{h_\mu\left(T,\mathscr{A}\right)\middle|\, \text{$\mathscr{A}$ is a finite measurable partition of $X$}\right\}. \]

We will need the following lemma later. 
See Theorem 4.12 (iv) of the book \cite[\S 4.5]{Walters book} for the proof.

\begin{lemma}  \label{lemma: Kolmogorov--Sinai entropy and conditional entropy}
If $\mathscr{A}$ and $\mathscr{A}^\prime$ are two finite measurable partitions of $X$, then
\[  h_\mu(T, \mathscr{A})  \leq h_\mu\left(T, \mathscr{A}^\prime\right) 
    + H_\mu\left(\mathscr{A} \middle| \mathscr{A}^\prime\right). \]
\end{lemma}

\section{Proof of $P^w_{\mathrm{var}}\left(\pi, T, f\right) \leq  P^w\left(\pi, T, f\right)$} \label{section: topology bounds measure}

Let $\pi:(X, T)\to (Y, S)$ be a factor map between dynamical systems, and let 
$f:X\to \mathbb{R}$ be a continuous function.
The purpose of this section is to prove a half of the variational principle:

\begin{proposition} \label{prop: topology bounds measure}
For any $0\leq w\leq 1$ and $\mu\in \mathscr{M}^T(X)$
\[ w h_\mu(T) + (1-w)h_{\pi_*\mu}(S) + w\int_X f\, d\mu  \leq P^w\left(\pi, T, f\right). \]
Therefore $P^w_{\mathrm{var}}\left(\pi, T, f\right) \leq  P^w\left(\pi, T, f\right)$.
\end{proposition}

\begin{proof}
Set $\nu = \pi_*\mu$. This is an invariant probability measure on $Y$.
We will prove 
\begin{equation} \label{eq: topology bounds measure}
 w h_\mu(T) + (1-w)h_{\nu}(S) + w\int_X f\, d\mu  \leq P^w\left(\pi, T, f\right) + 1+ 2\log 2.
\end{equation}
If this is proved, then we will get the above statement by the standard \textit{amplification trick}.
Namely, for each natural number $m$, 
we apply (\ref{eq: topology bounds measure}) to $\pi:(X, T^m)\to (Y, S^m)$ with a continuous function 
$\mathbb{S}_m f:X\to \mathbb{R}$:
\[ w h_\mu(T^m) + (1-w)h_{\nu}(S^m) + w\int_X \mathbb{S}_m f\, d\mu  \leq P^w\left(\pi, T^m, \mathbb{S}_m f\right) + 1+ 2\log 2. \]
We have $h_\mu(T^m) = m h_\mu(T)$, $h_{\nu}(S^m) = m h_{\nu}(S)$, $\int_X \mathbb{S}_m f\, d\mu = m\int_X f\, d\mu$ and 
\[ P^w\left(\pi, T^m, \mathbb{S}_m f\right) = m P^w\left(\pi, T, f\right) \quad (\text{Lemma \ref{lemma: amplification}}). \]
Hence 
\[ w h_\mu(T) + (1-w)h_{\nu}(S) + w\int_X f\, d\mu  \leq P^w\left(\pi, T, f\right) + \frac{1+ 2\log 2}{m}. \]
Letting $m\to \infty$, we get the statement.
So it is enough to prove (\ref{eq: topology bounds measure}).

Let $\mathscr{A} = \{A_1, \dots, A_\alpha\}$ be a finite measurable partition of $Y$, and let 
$\mathscr{B}$ be a finite measurable partition of $X$.
We will prove that 
\begin{equation}   \label{eq: half of variational principle for partitions A and B}
   w h_\mu\left(T,\mathscr{B}\right) + (1-w) h_\nu\left(S,\mathscr{A}\right) + 
    w \int_X f\, d\mu  \leq P^w\left(\pi, T, f\right) + 1 + 2\log 2. 
\end{equation}    
For each $A_a$ in $\mathscr{A}$ ($1\leq a \leq \alpha$), we take a compact subset $C_a\subset A_a$ satisfying 
\begin{equation}   \label{eq: choice of C_a}
  \sum_{a=1}^\alpha \nu\left(A_a\setminus C_a\right) <   \frac{1}{\log \alpha}. 
\end{equation}  
We set $C_0 = Y\setminus \left(C_1\cup \dots \cup C_\alpha\right)$ and 
$\mathscr{C} = \{C_0, C_1, C_2, \dots, C_\alpha\}$.

\begin{claim}  \label{claim: partitions A and C}
$\mathscr{C}$ is a finite measurable partition of $Y$ satisfying 
\[  h_\nu\left(S, \mathscr{A}\right) < h_\nu\left(S,\mathscr{C}\right) + 1. \]
\end{claim}

\begin{proof}
From Lemma \ref{lemma: Kolmogorov--Sinai entropy and conditional entropy}
\[ h_\nu\left(S, \mathscr{A}\right) \leq  h_\nu\left(S,\mathscr{C}\right) + H_\nu\left(\mathscr{A}\middle| \mathscr{C}\right). \]
Since $C_a\subset A_a$ for $1\leq a\leq \alpha$
\[ H_\nu\left(\mathscr{A}\middle| \mathscr{C}\right)  = \nu(C_0) 
    \sum_{a=1}^\alpha\left(-\frac{\nu(A_a\cap C_0)}{\nu(C_0)} \log \frac{\nu(A_a\cap C_0)}{\nu(C_0)}\right)
    \leq \nu(C_0) \log \alpha. \]
The last term is smaller than one by (\ref{eq: choice of C_a}).
\end{proof}

We consider $\mathscr{B}\vee \pi^{-1}(\mathscr{C})$, which has the form 
\[ \mathscr{B}\vee \pi^{-1}(\mathscr{C}) = \left\{B_{ab}\middle|\, 0\leq a \leq \alpha, 1\leq b\leq \beta_a\right\}, \quad 
    \pi^{-1}(C_a) = \bigcup_{b=1}^{\beta_a} B_{ab} \quad (0\leq a\leq \alpha). \]
For each $B_{ab}$ $(0\leq a\leq \alpha, 1\leq b\leq \beta_a)$ we take a compact subset $D_{ab}\subset B_{ab}$
such that
\begin{equation} \label{eq: choice of D_{ab}}
    \sum_{a=0}^\alpha \log \beta_a \left(\sum_{b=1}^{\beta_a} \mu\left(B_{ab}\setminus D_{ab}\right)\right) < 1.
\end{equation}
We set 
\[ D_{a0} = \pi^{-1}(C_a) \setminus \bigcup_{b=1}^{\beta_a} D_{ab}, \quad (0\leq a \leq \alpha). \]
We define 
\[  \mathscr{D} = \left\{D_{ab}\middle|\, 0\leq a\leq \alpha, 0\leq b\leq \beta_a \right\}. \]

\begin{claim}  \label{claim: partitions B and D}
$\mathscr{D}$ is a finite measurable partition of $X$ with $\pi^{-1}(\mathscr{C}) \prec \mathscr{D}$ and 
\[  h_\mu\left(T,\mathscr{B}\right) \leq h_\mu\left(T, \mathscr{D}\right) + 1. \]
\end{claim}

\begin{proof}
$\pi^{-1}(\mathscr{C}) \prec \mathscr{D}$ is obvious by the construction.
\begin{equation*}
   \begin{split}
        h_\mu\left(T,\mathscr{B}\right)  & \leq h_\mu\left(T, \mathscr{B}\vee \pi^{-1}(\mathscr{C})\right) \\
        &\leq h_\mu\left(T, \mathscr{D}\right) + H_\mu\left(\mathscr{B}\vee \pi^{-1}(\mathscr{C})\middle| \mathscr{D}\right)
        \quad  \text{by Lemma \ref{lemma: Kolmogorov--Sinai entropy and conditional entropy}}.
   \end{split}
\end{equation*}   
Since $D_{ab}\subset B_{ab}$ for $0\leq a\leq \alpha$ and $1\leq b\leq \beta_a$
\begin{equation*}
   \begin{split}
       H_\mu\left(\mathscr{B}\vee \pi^{-1}(\mathscr{C})\middle| \mathscr{D}\right)
     & = \sum_{a=0}^\alpha \mu(D_{a0}) \sum_{b=1}^{\beta_a}\left(-\frac{\mu(D_{a0}\cap B_{ab})}{\mu(D_{a0})}
       \log \frac{\mu(D_{a0}\cap B_{ab})}{\mu(D_{a0})}\right)  \\
     &\leq \sum_{a=0}^\alpha \mu(D_{a0}) \log \beta_a  \\
     & < 1 \quad \text{by (\ref{eq: choice of D_{ab}})}.
   \end{split}
\end{equation*}       
\end{proof}

We will prove that 
\[ w h_\mu\left(T, \mathscr{D}\right) + (1-w) h_\nu\left(S,\mathscr{C}\right) + 
    w \int_X f\, d\mu \leq P^w\left(\pi, T, f\right) + 2\log 2. \]
If this is proved, then (\ref{eq: half of variational principle for partitions A and B}) will follow from 
Claims \ref{claim: partitions A and C} and \ref{claim: partitions B and D}.

From the definition of the entropy,
\begin{equation*}
  \begin{split}
  w h_\mu\left(T, \mathscr{D}\right) + (1-w) h_\nu\left(S,\mathscr{C}\right) 
   & = \lim_{N\to \infty} \left(w \cdot \frac{H_\mu\left(\mathscr{D}^N\right)}{N} + (1-w)\cdot \frac{H_\nu\left(\mathscr{C}^N\right)}{N}\right) \\
   & = \lim_{N\to \infty}\frac{1}{N} 
   \left\{H_\nu\left(\mathscr{C}^N\right)  + w \left(H_\mu\left(\mathscr{D}^N\right) - H_\nu\left(\mathscr{C}^N\right)\right)\right\}.
  \end{split} 
\end{equation*}
Since $\nu = \pi_*\mu$, we have $H_\nu\left(\mathscr{C}^N\right) = H_\mu\left(\pi^{-1}\left(\mathscr{C}^N\right)\right)$.
Since $\pi^{-1}\left(\mathscr{C}^N\right) \prec \mathscr{D}^N$,
\[ H_\mu\left(\mathscr{D}^N\right) - H_\mu\left(\pi^{-1}\left(\mathscr{C}^N\right)\right) 
    = H_\mu\left(\mathscr{D}^N\middle| \pi^{-1}\left(\mathscr{C}^N\right)\right). \]
So 
\[ w h_\mu\left(T, \mathscr{D}\right) + (1-w) h_\nu\left(S,\mathscr{C}\right) = 
    \lim_{N\to \infty}\frac{1}{N} 
   \left\{H_\nu\left(\mathscr{C}^N\right) + w \cdot H_\mu\left(\mathscr{D}^N\middle| \pi^{-1}\left(\mathscr{C}^N\right)\right)\right\}. \]
We have 
\[ \int_X f\, d\mu = \frac{1}{N} \int_X \mathbb{S}_N f\, d\mu. \]
Therefore 
\begin{equation}  \label{eq: formula by conditional entropy}
  \begin{split}
   & w h_\mu\left(T, \mathscr{D}\right) + (1-w) h_\nu\left(S,\mathscr{C}\right) + w\int_X f\, d\mu  \\
   &  =  \lim_{N\to \infty}\frac{1}{N} 
   \left\{H_\nu\left(\mathscr{C}^N\right) + w \cdot H_\mu\left(\mathscr{D}^N\middle| \pi^{-1}\left(\mathscr{C}^N\right)\right)
   + w \int_X \mathbb{S}_N f\, d\mu \right\}
  \end{split} 
\end{equation}

For $C\in \mathscr{C}^N$, we define
\[ \mathscr{D}^N_C = \left\{D\in \mathscr{D}^N\middle| D\cap \pi^{-1}(C) \neq \emptyset\right\}
   = \left\{D\in \mathscr{D}^N\middle| D\subset\pi^{-1}(C) \right\}. \]
Then  
\[  \pi^{-1}(C) = \bigcup_{D\in \mathscr{D}^N_C} D. \]  
   
For $C\in \mathscr{C}^N$ with $\nu(C)>0$ and $D\in \mathscr{D}^N_C$, we set 
\[ \mu(D|C) = \frac{\mu(D)}{\nu(C)} = \frac{\mu(D)}{\mu\left(\pi^{-1}(C)\right)}. \]
For each $C\in \mathscr{C}^N$ with $\nu(C)>0$
we have 
\[  \sum_{D\in \mathscr{D}^N_C} \mu(D|C) = 1. \]

\begin{claim} \label{eq: bridge between measure and topology}
\begin{equation*}
   H_\nu\left(\mathscr{C}^N\right) + w \cdot H_\mu\left(\mathscr{D}^N\middle| \pi^{-1}\left(\mathscr{C}^N\right)\right)
   + w \int_X \mathbb{S}_N f\, d\mu  
   \leq \log \sum_{C\in \mathscr{C}^N} \left(\sum_{D\in \mathscr{D}^N_C} e^{\sup_{D} \mathbb{S}_N f}\right)^w.
\end{equation*}   
\end{claim}

\begin{proof}
We have 
\begin{equation*}
  \begin{split}
   \int_X \mathbb{S}_N f\, d\mu & = \sum_{D\in \mathscr{D}^N} \int_D \mathbb{S}_N f\, d\mu  
    \leq \sum_{D\in \mathscr{D}^N} \mu(D) \sup_{D} \mathbb{S}_N f  \\
    & = \sum_{\substack{C\in \mathscr{C}^N \\ \text{with $\nu(C)>0$}}} \nu(C) 
    \left(\sum_{D\in \mathscr{D}^N_C} \mu(D|C) \sup_D \mathbb{S}_N f\right).
  \end{split}
\end{equation*}  
Hence 
\begin{equation*}
  \begin{split}
  & H_\mu\left(\mathscr{D}^N\middle| \pi^{-1}\left(\mathscr{C}^N\right)\right) + \int_X \mathbb{S}_N f\, d\mu \\
  & \leq \sum_{\substack{C\in \mathscr{C}^N \\ \text{with $\nu(C)>0$}}} 
  \nu(C)\left\{ \sum_{D\in \mathscr{D}^N_C} \left( -\mu(D|C)\log \mu(D|C)
     + \mu(D|C) \sup_D \mathbb{S}_N f\right)\right\} 
  \end{split}
\end{equation*}  
By Lemma \ref{lemma: calculus} (2) 
\[ \sum_{D\in \mathscr{D}^N_C} \left( -\mu(D|C)\log \mu(D|C)
     + \mu(D|C) \sup_D \mathbb{S}_N f\right) \leq 
    \log  \sum_{D\in \mathscr{D}^N_C} e^{\sup_D \mathbb{S}_N f}. \]
So
\[ H_\mu\left(\mathscr{D}^N\middle| \pi^{-1}\left(\mathscr{C}^N\right)\right) + \int_X \mathbb{S}_N f\, d\mu 
\leq \sum_{C\in \mathscr{C}^N} \nu(C) 
 \left(\log \sum_{D\in \mathscr{D}^N_C} e^{\sup_D \mathbb{S}_N f}\right). \]
Therefore
\begin{equation*}
  \begin{split}
    & H_\nu\left(\mathscr{C}^N\right) + w \cdot H_\mu\left(\mathscr{D}^N\middle| \pi^{-1}\left(\mathscr{C}^N\right)\right)
   + w \int_X \mathbb{S}_N f\, d\mu     \\
   & \leq    \sum_{C\in \mathscr{C}^N} \left\{-\nu(C) \log \nu(C) + \nu(C) 
      \log \left(\sum_{D\in \mathscr{D}^N_C} e^{\sup_D \mathbb{S}_N f}\right)^w\right\} \\
   & \leq    \log \sum_{C\in \mathscr{C}^N} \left(\sum_{D\in \mathscr{D}^N_C} e^{\sup_{D} \mathbb{S}_N f}\right)^w 
   \quad \text{by Lemma \ref{lemma: calculus} (2)}.
  \end{split} 
\end{equation*}
\end{proof}

We take metrics $d$ and $d^\prime$ on $X$ and $Y$ respectively.
Recall that $C_a$ $(1\leq a\leq \alpha)$ are mutually disjoint compact subsets of $Y$ and that 
$D_{ab}$ $(0\leq a \leq \alpha, 1\leq b \leq \beta_a)$ are mutually disjoint compact subsets of $X$.
Hence we can take $\varepsilon>0$ such that
\begin{itemize}
 \item for any $y\in C_a$ and $y^\prime\in C_{a^\prime}$ with distinct $1\leq a, a^\prime\leq \alpha$
  \[  \varepsilon < d^\prime(y, y^\prime). \]
 \item for any $x\in D_{ab}$ and $x^\prime\in D_{a b^\prime}$ with $0\leq a\leq \alpha$ and distinct $1\leq b, b^\prime\leq \beta_a$
  \[  \varepsilon < d(x, x^\prime). \]
\end{itemize}

\begin{claim}  \label{claim: bounds on overlap}
Let $N$ be a natural number.
  \begin{enumerate}
   \item[(1)]  If a subset $V\subset Y$ has $\diam(V,d^\prime_N) < \varepsilon$ then the number of $C\in \mathscr{C}^N$
   having nonempty intersection with $V$ is at most $2^N$:
   \[  \left|\{C\in \mathscr{C}^N| C\cap V\neq \emptyset\}\right| \leq 2^N. \]
   \item[(2)]  If a subset $U\subset X$ has $\diam(U, d_N) < \varepsilon$ then for each $C\in \mathscr{C}^N$ the number of
   $D\in \mathscr{D}^N_C$ having nonempty intersection with $U$ is at most $2^N$:
   \[  \left|\{D\in \mathscr{D}^N_C| D\cap U \neq \emptyset\}\right| \leq 2^N. \]
  \end{enumerate} 
\end{claim}

\begin{proof}
(1) For each $0\leq k < N$, the set $S^k V$ may have nonempty intersection with $C_0$ and at most one 
set in $\{C_1, C_2, \dots, C_\alpha\}$. The above statement follows from this.

(2) Suppose $C\in \mathscr{C}^N$ has the form 
\[   C = C_{a_0}\cap S^{-1}C_{a_1}\cap S^{-2}C_{a_2} \cap \dots \cap S^{-(N-1)}C_{a_{N-1}}  \]
with $0\leq a_0, \dots, a_{N-1} \leq \alpha$.
Recall that $\{D_{a_k 0}, D_{a_k 1}, D_{a_k 2}, \dots, D_{a_k \beta_{a_k}}\}$ is a partition of $\pi^{-1}(C_{a_k})$.
Then any set $D\in \mathscr{D}^N_C$ has the form 
\[ D = D_{a_0 b_0} \cap T^{-1}D_{a_1 b_1} \cap T^{-2} D_{a_2 b_2} \cap \dots \cap T^{-(N-1)} D_{a_{N-1} b_{N-1}} \]
with $0 \leq b_k \leq \beta_{a_k}$ for $0\leq k \leq N-1$. 

For each $0\leq k <N$, the set $T^k U$ may have nonempty intersection with $D_{a_k 0}$ and at most one set in 
$\{D_{a_k 1}, D_{a_k 2}, \dots, D_{a_k \beta_{a_k}}\}$. 
Now the above statement follows from this.
\end{proof}

Let $N$ be a natural number.
Suppose we are given an open cover $Y= V_1\cup \dots \cup V_n$ with $\diam\left(V_i, d^\prime_N\right)  < \varepsilon$
for all $1\leq i\leq n$.
Moreover, suppose that for each $1\leq i\leq n$ we are given an open cover 
$\pi^{-1}(V_i) = U_{i1}\cup U_{i2}\cup \dots \cup U_{i m_i}$ with 
$\diam\left(U_{ij}, d_N\right) < \varepsilon$ for all $1\leq j\leq m_i$.
We are going to prove 
\begin{equation}  \label{eq: bound by weighted pressure for U and V}
    \log \sum_{C\in \mathscr{C}^N} \left(\sum_{D\in \mathscr{D}^N_C} e^{\sup_{D} \mathbb{S}_N f}\right)^w  \leq 
   2N\log 2 + \log\sum_{i=1}^n\left(\sum_{j=1}^{m_i} e^{\sup_{U_{ij}} \mathbb{S}_N f}\right)^w.
\end{equation}  
Suppose this is proved. 
Then by Claim \ref{eq: bridge between measure and topology}
\[ H_\nu\left(\mathscr{C}^N\right) + w \cdot H_\mu\left(\mathscr{D}^N\middle| \pi^{-1}\left(\mathscr{C}^N\right)\right)
   + w \int_X \mathbb{S}_N f\, d\mu 
   \leq 2N\log 2 + \log\sum_{i=1}^n\left(\sum_{j=1}^{m_i} e^{\sup_{U_{ij}} \mathbb{S}_N f}\right)^w. \]
Taking the infimum over $\{V_i\}$ and $\{U_{ij}\}$ satisfying the above assumptions, we have 
\[ H_\nu\left(\mathscr{C}^N\right) + w \cdot H_\mu\left(\mathscr{D}^N\middle| \pi^{-1}\left(\mathscr{C}^N\right)\right)
   + w \int_X \mathbb{S}_N f\, d\mu 
   \leq 2N\log 2 + \log P^w\left(\pi, f, N, \varepsilon\right). \]
Divide this by $N$ and let $N\to \infty$. 
Recalling (\ref{eq: formula by conditional entropy}), we get 
\[ w h_\mu\left(T,\mathscr{D}\right) + (1-w) h_\nu\left(S,\mathscr{C}\right) +  w\int_X f\, d\mu 
   \leq 2 \log 2 + \lim_{N\to \infty} \frac{\log P^w\left(\pi, f, N, \varepsilon\right)}{N}. \]
Letting $\varepsilon\to 0$, we get the desired result:
\[ w h_\mu\left(T,\mathscr{D}\right) + (1-w) h_\nu\left(S,\mathscr{C}\right) +  w\int_X f\, d\mu 
   \leq 2 \log 2 + P^w\left(\pi, T, f\right). \] 
So the rest of the work is to prove (\ref{eq: bound by weighted pressure for U and V}).

For $D\in \mathscr{D}^N$, we have
\[ e^{\sup_D \mathbb{S}_N f} \leq \sum_{U_{ij}\cap D \neq \emptyset} e^{\sup_{U_{ij}} \mathbb{S}_N f}. \]
Here the sum is taken over the index $(i, j)$ such that $U_{ij}$ has nonempty intersection with $D$.

Let $C\in \mathscr{C}^N$. We define $\mathscr{V}_C$ as the set of $1\leq i \leq n$ such that 
$V_i \cap C \neq \emptyset$.
By Claim \ref{claim: bounds on overlap} (2), 
\[  \sum_{D\in \mathscr{D}_C^N} e^{\sup_{D} \mathbb{S}_N f} 
  \leq 2^N \sum_{i\in \mathscr{V}_C} \sum_{j=1}^{m_i} e^{\sup_{U_{ij}} \mathbb{S}_N f}. \]
Then (recall $0\leq w\leq 1$)
\begin{equation*}
  \begin{split}
    \left(\sum_{D\in \mathscr{D}_C^N} e^{\sup_{D} \mathbb{S}_N f}\right)^w  & \leq
    2^{Nw} \left(\sum_{i\in \mathscr{V}_C}\sum_{j=1}^{m_i} e^{\sup_{U_{ij}} \mathbb{S}_N f}\right)^w \\
    & \leq  2^{Nw} \sum_{i\in \mathscr{V}_C} \left(\sum_{j=1}^{m_i} e^{\sup_{U_{ij}} \mathbb{S}_N f}\right)^w \quad 
    \text{by Lemma \ref{lemma: calculus} (1)}.
  \end{split}
\end{equation*}
Hence
\[ \sum_{C\in \mathscr{C}^N} \left(\sum_{D\in \mathscr{D}_C^N} e^{\sup_{D} \mathbb{S}_N f}\right)^w \leq 
    2^{Nw} \sum_{C\in \mathscr{C}^N}
    \left\{\sum_{i\in \mathscr{V}_C} \left(\sum_{j=1}^{m_i} e^{\sup_{U_{ij}} \mathbb{S}_N f}\right)^w\right\}. \]
By Claim \ref{claim: bounds on overlap} (1), for each $1\leq i \leq n$, 
the number of $C\in\mathscr{C}^N$ satisfying $i \in \mathscr{V}_C$
is at most $2^N$. So the right-hand side is bounded from above by 
\[ 2^{Nw} \cdot 2^N \sum_{i=1}^n \left(\sum_{j=1}^{m_i} e^{\sup_{U_{ij}} \mathbb{S}_N f}\right)^w. \]
Therefore
\[ \sum_{C\in \mathscr{C}^N} \left(\sum_{D\in \mathscr{D}_C^N} e^{\sup_{D} \mathbb{S}_N f}\right)^w 
     \leq 2^{Nw} \cdot 2^N \sum_{i=1}^n \left(\sum_{j=1}^{m_i} e^{\sup_{U_{ij}} \mathbb{S}_N f}\right)^w. \]
Taking the logarithm,
\begin{equation*}
  \begin{split}
   \log \sum_{C\in \mathscr{C}^N} \left(\sum_{D\in \mathscr{D}_C^N} e^{\sup_{D} \mathbb{S}_N f}\right)^w 
    & \leq (N+Nw) \log 2 + \log \sum_{i=1}^n \left(\sum_{j=1}^{m_i} e^{\sup_{U_{ij}} \mathbb{S}_N f}\right)^w \\
    & \leq 2N\log 2 + \log \sum_{i=1}^n \left(\sum_{j=1}^{m_i} e^{\sup_{U_{ij}} \mathbb{S}_N f}\right)^w.
   \end{split}
\end{equation*} 
This is the estimate (\ref{eq: bound by weighted pressure for U and V}).
So we have finished the proof of the proposition.
\end{proof}

\section{Zero dimensional principal extension}  \label{section: zero dimensional principal extension}

In this section we prepare some definitions and results on principal extensions.
The main reference is the book of Downarowicz \cite{Downarowicz}.

Let $\pi:(X, T)\to (Y, S)$ be a factor map between dynamical systems.
Let $d$ be a metric on $X$.
We define the \textbf{topological conditional entropy} of $\pi$ by 
\[ h_{\mathrm{top}}(X, T|Y, S) 
= \lim_{\varepsilon\to 0} \left(\lim_{N\to \infty} \frac{\sup_{y\in Y} \log \#\left(\pi^{-1}(y), N, \varepsilon\right)}{N}\right). \]
Here $\#\left(\pi^{-1}(y), N, \varepsilon\right)$ is the number defined by (\ref{eq: covering number}).
It is easy to check that the quantity
\[ \sup_{y\in Y} \log \#\left(\pi^{-1}(y), N, \varepsilon\right) \]
 is sub-additive in $N$ and monotone in $\varepsilon$.
So the above limits exist.
This definition of the topological conditional entropy is due to \cite[Lemma 6.8.2]{Downarowicz}.

The factor map $\pi$ is said to be \textbf{principal} if $h_{\mathrm{top}}(X,T|Y,S)=0$.
In the case that this condition holds, 
the dynamical system $(X, T)$ is called a \textbf{principal extension} of $(Y, S)$.

The next theorem shows an important consequence of this condition.
This is proved in \cite[Corollary 6.8.9]{Downarowicz}. 
(See also the paper of Ledrappier--Walters \cite{Ledrappier--Walters}.)

\begin{theorem} \label{theorem: conditional variational principle and principal extension}
A principal factor map preserves Kolmogorov--Sinai entropy.
Namely, if $\pi:(X, T)\to (Y, S)$ is a principal factor map between dynamical systems then
for any invariant probability measure $\mu \in \mathscr{M}^T(X)$ 
\[  h_\mu(T) = h_{\pi_*\mu}(S). \]
\end{theorem}

\begin{remark}
Indeed, \cite[Corollary 6.8.9]{Downarowicz} proves the following more precise result:
Let $\pi:(X, T)\to (Y, S)$ be a factor map with $\h(Y, S) < \infty$.
Then $\pi$ is a principal factor map if and only if $h_\mu(T) = h_{\pi_*\mu}(S)$ for all $\mu\in \mathscr{M}^T(X)$.
\end{remark}

\begin{lemma}  \label{lemma: fiber product and principal extension}
Let $(X, T), (Y, S), (Y^\prime, S^\prime)$ be dynamical systems.
Let $\pi:X\to Y$ be a factor map, and let $\phi:Y^\prime \to Y$ be a principal factor map.
We define the fiber product (see Lemma \ref{lemma: fiber product})
\[ X\times_Y Y^\prime = \left\{(x, y)\in X\times Y^\prime\middle|\, \pi(x) = \phi(y)\right\}. \]
$\left(X\times_Y Y^\prime, T\times S^\prime\right)$ becomes a dynamical system.
We define factor maps $\varphi: X\times_Y Y^\prime\to X$ and $\Pi:X\times_Y Y^\prime\to Y^\prime$ by 
\[ \varphi(x, y) = x, \quad \Pi(x,y) = y. \]
\[  \xymatrix{ \left(X\times_Y Y^\prime, T\times S^\prime\right) \ar[rrr]^{\varphi} \ar[d]_{\Pi} & & & (X,T) \ar[d]^{\pi} \\
                       (Y^\prime, S^\prime) \ar[rrr]_{\phi: \text{ principal}} & & & (Y,S) } \]
Then $\varphi$ is a principal factor map. 
(The map $\Pi$ is not used in this statement, but we have introduced it for the convenience in the sequel.)
\end{lemma}

\begin{proof}
Let $d$ and $d^\prime$ be metrics on $X$ and $Y^\prime$ respectively.
We define a metric $\rho$ on $X\times_Y Y^\prime$ by 
\[ \rho\left((x_1, y_1), (x_2, y_2)\right) = 
    \max\left(d(x_1, x_2), d^\prime(y_1, y_2)\right). \]
For any natural number $N$ and $x\in X$, the metric space 
\[  \left(\varphi^{-1}(x), \rho_N\right) \]
is isometric to $\left(\phi^{-1}\left(\pi(x)\right), d^\prime_N\right)$.
Therefore for any $\varepsilon>0$
\[ \#\left(\varphi^{-1}(x), N, \varepsilon\right) = \#\left(\phi^{-1}\left(\pi(x)\right), N, \varepsilon\right). \]
So (recall that a factor map is always surjective)
\[ \sup_{x\in X} \#\left(\varphi^{-1}(x), N, \varepsilon\right) = \sup_{x\in X} \#\left(\phi^{-1}(\pi(x)), N, \varepsilon\right)
    = \sup_{y\in Y} \#\left(\phi^{-1}(y), N, \varepsilon\right). \]
Thus
\[ \h\left(X\times_Y Y^\prime, T\times S^\prime\middle| X, T\right) 
     = \h\left(Y^\prime, S^\prime\middle| Y, S\right) = 0. \]
\end{proof}

The next theorem is a key technical result.
This is proved in \cite[Theorem 7.6.1]{Downarowicz}.
(See also \cite{Downarowicz--Huczek}.)
Here recall that a compact metrizable space is said to be \textbf{zero dimensional}
if clopen subsets\footnote{A subset of a topological space is called \textbf{clopen} if it is closed and open.}
form an open basis of the topology.
For example, the Cantor set $\{0, 1\}^{\mathbb{N}}$ is zero dimensional.
A dynamical system $(X, T)$ is said to be zero dimensional if $X$ is a zero dimensional compact metrizable space.

\begin{theorem}\label{theorem: zero dimensional principal extension}
Every dynamical system has a zero dimensional principal extension.
Namely, for any dynamical system $(X, T)$, there exist a dynamical system $(X^\prime, T^\prime)$ and 
a factor map $\phi:X^\prime\to X$ such that $X^\prime$ is zero dimensional and $\phi$ is principal.
\end{theorem}

Recall that
we have defined two terms $P^w(\pi,T, f)$ and $P^w_{\mathrm{var}}(\pi, T, f)$
in \S \ref{section: weighted topological pressure}.

\begin{corollary}  \label{cor: reduction to zero dimensional case}
Let $\pi:(X, T)\to (Y, S)$ be a factor map between dynamical systems with a continuous function 
$f:X\to \mathbb{R}$.
There exists a factor map $\pi^\prime:(X^\prime, T^\prime)\to (Y^\prime, S^\prime)$ with a continuous function 
$f^\prime:X^\prime\to \mathbb{R}$ satisfying the following two conditions.
\begin{enumerate}
   \item[(1)] $X^\prime$ and $Y^\prime$ are zero dimensional.
   \item[(2)] For any $0\leq w\leq 1$ we have
   \[ P^w(\pi,T, f) \leq P^w(\pi^\prime, T^\prime, f^\prime), \quad 
       P^w_{\mathrm{var}}(\pi^\prime, T^\prime, f^\prime) \leq P^w_{\mathrm{var}}(\pi,T,  f). \]
\end{enumerate}
\end{corollary}

\begin{proof}
By Theorem \ref{theorem: zero dimensional principal extension}, there exists a zero dimensional 
principal extension $\phi:(Y^\prime,S^\prime)\to (Y, S)$.
We consider the fiber product $\left(X\times_Y Y^\prime, T\times S^\prime\right)$ and 
the projections $\varphi: X\times_Y Y^\prime\to X$ and $\Pi: X\times_Y Y^\prime\to Y^\prime$
as in Lemma \ref{lemma: fiber product and principal extension}.
Then $\varphi$ is a principal factor map.
\[  \xymatrix{ \left(X\times_Y Y^\prime, T\times S^\prime\right) \ar[rrr]^{\varphi} \ar[d]_{\Pi} & & & (X,T) \ar[d]^{\pi} \\
                       (Y^\prime, S^\prime) \ar[rrr]_{\phi: \text{ principal}} & & & (Y,S) } \]

By Lemma \ref{lemma: fiber product}, for any $0\leq w\leq 1$
\[  P^w(\pi, T, f) \leq P^w\left(\Pi, T\times S^\prime, f\circ\varphi\right). \]
Here the right-hand side is $P^w\left(\Pi, X\times_Y Y^\prime, T\times S^\prime, Y^\prime, S^\prime, f\circ \varphi\right)$.
By Theorem \ref{theorem: conditional variational principle and principal extension}, for any invariant probability measure
$\mu \in \mathscr{M}^{T\times S^\prime}\left(X\times_Y Y^\prime\right)$
\[   h_\mu(T\times S^\prime) = h_{\varphi_*\mu}(T), \quad 
      h_{\Pi_*\mu}(S^\prime) = h_{\phi_*\Pi_*\mu}(S) = h_{\pi_*\varphi_*\mu}(S).       \] 
Then\footnote{Here we prove $P^w_{\mathrm{var}}\left(\Pi, T\times S^\prime, f\circ\varphi\right)
\leq P^w_{\mathrm{var}}(\pi, T, f)$.
Indeed we can prove the equality
$P^w_{\mathrm{var}}\left(\Pi, T\times S^\prime, f\circ\varphi\right)
= P^w_{\mathrm{var}}(\pi, T, f)$
because the map $\varphi_*:\mathscr{M}^{T\times S^\prime}(X\times_Y Y^\prime)\to \mathscr{M}^{T}(X)$ is surjective. 
But we do not need this.}
\begin{equation}  \label{eq: fiber product and P_var}
  \begin{split}
    & P^w_{\mathrm{var}}\left(\Pi, T\times S^\prime, f\circ\varphi\right)  \\
   & =\sup_{\mu\in \mathscr{M}^{T\times S^\prime}(X\times_Y Y^\prime)}
       \left\{w h_\mu(T\times S^\prime) + (1-w)h_{\Pi_*\mu}(S^\prime) + w\int_{X\times_Y Y^\prime} 
       f\circ \varphi\, d\mu\right\} \\
    & =  \sup_{\mu\in \mathscr{M}^{T\times S^\prime}(X\times_Y Y^\prime)}
       \left\{w h_{\varphi_*\mu}(T) + (1-w) h_{\pi_*\varphi_*\mu}(S) + w \int_X f \, d\left(\varphi_*\mu\right)\right\}  \\
    & \leq P^w_{\mathrm{var}}(\pi, T, f)
  \end{split}  
\end{equation}

By applying Theorem \ref{theorem: zero dimensional principal extension} to the system $\left(X\times_Y Y^\prime, T\times S^\prime\right)$,
there exists a zero dimensional principal extension 
$\psi: (X^\prime, T^\prime) \to \left(X\times_Y Y^\prime, T\times S^\prime\right)$.
\[  \xymatrix{ (X^\prime, T^\prime) \ar[rrr]^{\psi: \text{ principal}} \ar[rrrd]_{\Pi\circ \psi} & & &
    \left(X\times_Y Y^\prime, T\times S^\prime\right) \ar[rrr]^{\varphi: \text{ principal}} \ar[d]_{\Pi} & & & (X,T) \ar[d]^{\pi} \\
         &  &    &        (Y^\prime, S^\prime) \ar[rrr]_{\phi: \text{ principal}} & & & (Y,S) } \]

By Lemma \ref{lemma: pressure and factor map}
\[    P^w\left(\Pi, T\times S^\prime, f\circ\varphi\right)    \leq  P^w\left(\Pi\circ \psi, T^\prime, f\circ \varphi\circ \psi\right). \]
Here the right-hand side is 
$P^w\left(\Pi\circ \psi, X^\prime, T^\prime, Y^\prime, S^\prime, f \circ \varphi\circ \psi\right)$.
As in the above (\ref{eq: fiber product and P_var}), by Theorem \ref{theorem: conditional variational principle and principal extension},
\[   P^w_{\mathrm{var}}\left(\Pi\circ \psi, T^\prime, f\circ \varphi\circ \psi\right) \leq P^w\left(\Pi, T\times S^\prime, f\circ\varphi\right). \]

So we conclude 
\[ P^w(\pi, T, f) \leq   P^w\left(\Pi\circ \psi, T^\prime, f\circ \varphi\circ \psi\right), \quad 
    P^w_{\mathrm{var}}\left(\Pi\circ \psi, T^\prime, f\circ \varphi\circ \psi\right)  \leq  P^w_{\mathrm{var}}(\pi, T, f). \]
Set $\pi^\prime := \Pi\circ \psi: (X^\prime, T^\prime) \to (Y^\prime, S^\prime)$ and 
$f^\prime := f\circ \varphi\circ \psi: X^\prime \to \mathbb{R}$.
These satisfy the required conditions.
\end{proof}

\section{Completion of the proof of the variational principle}  \label{section: completion of the proof of the variational principle}

In this section we prove $P^w(\pi, T, f) \leq P^w_{\mathrm{var}}(\pi, T, f)$ and complete the proof of the variational principle.
First we consider the case of zero dimensional dynamical systems.
Later we will reduce the general case to this zero dimensional case.

\begin{proposition} \label{prop: main result in zero dimensional case}
Let $\pi:(X,T)\to (Y, S)$ be a factor map between \textbf{zero dimensional} dynamical systems.
Then for any $0\leq w\leq 1$ and a continuous function $f:X\to \mathbb{R}$
\[  P^w(\pi, T, f) \leq P^w_{\mathrm{var}}(\pi, T, f). \]
\end{proposition}

\begin{proof}
Let $\varepsilon>0$.
We will prove that there exists $\mu\in \mathscr{M}^T(X)$ satisfying 
\[ w h_\mu(T) + (1-w)h_{\pi_*\mu}(S) + w \int_X f\, d\mu 
   \geq \lim_{N\to \infty} \frac{\log P^w\left(\pi, f, N, \varepsilon\right)}{N}. \]

We take metrics $d$ and $d^\prime$ on $X$ and $Y$ respectively.
Let $Y = A_1\cup \dots \cup A_\alpha$ be a clopen partition (i.e. $A_a$ are mutually disjoint clopen subsets of $Y$)
with $\diam \left(A_a, d^\prime\right) < \varepsilon$ for all $1\leq a\leq \alpha$. Here we have used $\dim Y=0$.

From $\dim X = 0$,
for each $1\leq a\leq \alpha$, we can also take a clopen partition 
\[ \pi^{-1}(A_a) = \bigcup_{b=1}^{\beta_a} B_{ab} \quad 
    \text{with $\diam\left(B_{ab}, d\right) < \varepsilon$ for all $1\leq b\leq \beta_a$}. \]
Set $\mathscr{A} = \{A_1, \dots, A_\alpha\}$ and 
$\mathscr{B} = \left\{B_{ab}\middle|\, 1\leq a\leq \alpha, 1\leq b\leq \beta_a\right\}$.
These are clopen partitions of $Y$ and $X$ respectively.
We have $\mathscr{B} \prec \pi^{-1}(\mathscr{A})$.

Let $N$ be a natural number. We have $\mathscr{B}^N \prec \pi^{-1}\left(\mathscr{A}^N\right)$.
For each nonempty $A\in \mathscr{A}^N$ we define 
\[ \mathscr{B}^N_A = \left\{B\in \mathscr{B}^N\middle|\, B\cap \pi^{-1}(A) \neq \emptyset\right\}
    = \left\{B\in \mathscr{B}^N\middle|\, B\subset \pi^{-1}(A), B \neq \emptyset\right\}. \]    
We have 
\[   \pi^{-1}(A) = \bigcup_{B\in \mathscr{B}^N_A} B. \]    
We set 
\[ Z_{N,A} = \sum_{B\in \mathscr{B}^N_A} e^{\sup_{B} \mathbb{S}_N f}. \]
Define
\[  Z_N = \sum_{A\in \mathscr{A}^N} \left(Z_{N,A}\right)^w. \]    
Here the sum is taken over only nonempty $A \in \mathscr{A}^N$.
When we consider below a sum over $A\in \mathscr{A}^N$ (or $B\in \mathscr{B}^N$), we always assume that 
$A$ (or $B$) is not empty.

We have 
\[  P^w\left(\pi, f, N, \varepsilon\right) \leq Z_N. \]
So it is enough to prove that there exists $\mu\in \mathscr{M}^T(X)$ satisfying 
\[ w h_\mu\left(T,\mathscr{B}\right) + (1-w)h_{\pi_*\mu}\left(S,\mathscr{A}\right) + w \int_X f\, d\mu 
   \geq \lim_{N\to \infty} \frac{\log Z_N}{N}, \]
where the limit in the right-hand side exists because $Z_N$ is sub-multiplicative in $N$.   
    
Let $N$ be a natural  number. 
For nonempty $B\in \mathscr{B}^N$, we denote by $\mathscr{A}^N(B)$ the unique element of $\mathscr{A}^N$
containing $\pi(B)$.  
For nonempty $A\in \mathscr{A}^N$, we have $\mathscr{A}^N(B) = A$ for all $B\in \mathscr{B}^n_A$.

For each nonempty set $B$ in $\mathscr{B}^N$ we take a point $x_B\in B$ satisfying 
$\mathbb{S}_N f(x_B) = \sup_{B} \mathbb{S}_N f$.
(Such a point exists because $B$ is closed.)
We define a probability measure on $X$ by 
\begin{equation*}
   \begin{split}
    \sigma_N & = \frac{1}{Z_N} \sum_{B\in \mathscr{B}^N} \left(Z_{N,\mathscr{A}^N(B)}\right)^{w-1} 
    e^{\mathbb{S}_N f(x_B)} \cdot \delta_{x_B} \\
    & = \frac{1}{Z_N} \sum_{A\in \mathscr{A}^N} \sum_{B\in \mathscr{B}^N_A} \left(Z_{N,A}\right)^{w-1} e^{\mathbb{S}_N f(x_B)} \cdot 
         \delta_{x_B}.
   \end{split}
\end{equation*}    
Here $\delta_{x_B}$ is the delta probability measure at the point $x_B$.
$\sigma_N$ is not an invariant measure in general.
We set 
\[ \mu_N = \frac{1}{N} \sum_{n=0}^{N-1} T^n_*\sigma_N. \] 
We can take a subsequence $\{\mu_{N_k}\}$ converging to an invariant probability measure $\mu$ on $X$
in the weak$^*$ topology.
We will prove that this measure $\mu$ satisfies
\[ w h_\mu(T) + (1-w)h_{\pi_*\mu}(S) + w \int_X f\, d\mu 
   \geq \lim_{N\to \infty} \frac{\log Z_N}{N}. \]

\begin{claim}  \label{claim: Z_n and measure sigma_N}
For every natural number $N$
\[  w H_{\sigma_N}\left(\mathscr{B}^N\right) + (1-w) H_{\pi_* \sigma_N}\left(\mathscr{A}^N\right) 
     + w \int_X \mathbb{S}_N f\, d\sigma_N = \log Z_N. \]
\end{claim}    
    
\begin{proof}
We have 
\[ \pi_*\sigma_N = \frac{1}{Z_N} \sum_{B\in \mathscr{B}^N} \left(Z_{N,\mathscr{A}^N(B)}\right)^{w-1} e^{\mathbb{S}_N f(x_B)}\cdot 
                           \delta_{\pi(x_B)}. \]
For each nonempty $A\in \mathscr{A}^N$
\begin{equation*}
  \begin{split}
  \pi_*\sigma_N(A) & = \frac{1}{Z_N} \sum_{B\in \mathscr{B}^N_A} \left(Z_{N,\mathscr{A}^N(B)}\right)^{w-1} e^{\mathbb{S}_N f(x_B)} \\
      &= \frac{1}{Z_N} \left(Z_{N,A}\right)^w \quad \text{by $\mathscr{A}^N(B) = A$ for $B\in \mathscr{B}^N_A$}.
  \end{split}
\end{equation*}  
Then 
\begin{equation}  \label{eq: entropy of pi_sigma_N for A^N}
  H_{\pi_*\sigma_N} \left(\mathscr{A}^N\right) = \log Z_N - w\sum_{A\in \mathscr{A}^N}\frac{\left(Z_{N, A}\right)^w}{Z_N}
     \log Z_{N,A}.
\end{equation}

For nonempty $B\in \mathscr{B}^N$
\[ \sigma_N(B) = \frac{\left(Z_{N, \mathscr{A}^N(B)}\right)^{w-1}}{Z_N} e^{\mathbb{S}_N f(x_B)}. \]
Then 
\begin{equation*}
  \begin{split}
   H_{\sigma_N}\left(\mathscr{B}^N\right) 
   =& - \sum_{B\in \mathscr{B}^N} \frac{\left(Z_{N, \mathscr{A}^N(B)}\right)^{w-1}}{Z_N} e^{\mathbb{S}_N f(x_B)} 
    \log \left(\frac{\left(Z_{N, \mathscr{A}^N(B)}\right)^{w-1}}{Z_N} e^{\mathbb{S}_N f(x_B)}\right) \\
    = &   \frac{\log Z_N}{Z_N} 
   \underbrace{\sum_{B\in \mathscr{B}^N} \left(Z_{N, \mathscr{A}^N(B)}\right)^{w-1} e^{\mathbb{S}_N f(x_B)}}_{(I)} \\
    & - \frac{w-1}{Z_N} \underbrace{\sum_{B\in \mathscr{B}^N} \left(Z_{N,\mathscr{A}^N(B)}\right)^{w-1} e^{\mathbb{S}_N f(x_B)}
      \log Z_{N, \mathscr{A}^N(B)}}_{(II)} \\
   &- \underbrace{\sum_{B\in \mathscr{B}^N} 
   \frac{\left(Z_{N, \mathscr{A}^N(B)}\right)^{w-1}}{Z_N} e^{\mathbb{S}_N f(x_B)} \mathbb{S}_N f (x_B)}_{(III)}.
  \end{split}
\end{equation*}  
We calculate the term $(I)$ by 
\[ (I) = \sum_{A\in \mathscr{A}^N} \sum_{B\in \mathscr{B}^N_A} \left(Z_{N, A}\right)^{w-1} e^{\mathbb{S}_N f (x_B)}
       =  \sum_{A\in \mathscr{A}^N} \left(Z_{N, A}\right)^{w-1} \cdot Z_{N,A}
       = Z_N. \]
The term $(II)$ is calculated by
\begin{equation*}
   (II) = \sum_{A\in \mathscr{A}^N} \sum_{B\in \mathscr{B}_A^N} \left(Z_{N,A}\right)^{w-1} e^{\mathbb{S}_N f(x_B)}
           \log Z_{N,A} 
        = \sum_{A\in \mathscr{A}^N} \left(Z_{N,A}\right)^w \log Z_{N,A}.
\end{equation*}
For the term $(III)$ we consider 
\[ \int_X \mathbb{S}_N f\, d\sigma_N  = \frac{1}{Z_N} \sum_{B\in \mathscr{B}^N} \left(Z_{N, \mathscr{A}^N(B)}\right)^{w-1}
                e^{\mathbb{S}_N f(x_B)} \mathbb{S}_N f(x_B) = (III). \]
Thus 
\[ H_{\sigma_N}\left(\mathscr{B}^N\right) + \int_X \mathbb{S}_N f\, d\sigma_N
    = \log Z_N - \frac{w-1}{Z_N} \sum_{A\in \mathscr{A}^N} \left(Z_{N,A}\right)^{w} \log Z_{N, A}. \]
Combining this with (\ref{eq: entropy of pi_sigma_N for A^N}) we get
\[ w H_{\sigma_N}\left(\mathscr{B}^N\right) + (1-w) H_{\pi_* \sigma_N}\left(\mathscr{A}^N\right) 
     + w \int_X \mathbb{S}_N f\, d\sigma_N = \log Z_N. \]
\end{proof}

\begin{claim}  \label{claim: entropy bounded from below}
Let $M$ and $N$ be natural numbers. We have
  \begin{equation*}
    \begin{split}
      \frac{1}{M} H_{\mu_N}\left(\mathscr{B}^M\right)  & \geq \frac{1}{N} H_{\sigma_N}\left(\mathscr{B}^N\right) 
       -\frac{2M \log |\mathscr{B}|}{N}, \\
      \frac{1}{M} H_{\pi_* \mu_N}\left(\mathscr{A}^M\right)  & \geq \frac{1}{N} H_{\pi_*\sigma_N}\left(\mathscr{A}^N\right) 
       -\frac{2M \log |\mathscr{A}|}{N}.
    \end{split}
  \end{equation*}
  Here $|\mathscr{A}|$ and $|\mathscr{B}|$ are the cardinalities of $\mathscr{A}$ and $\mathscr{B}$ respectively.
\end{claim}

\begin{proof}
This is rather standard. (See the proof of the standard variational principle in \cite[\S 8.2]{Walters book}.)
Here we provide the proof for $\mathscr{B}^M$. The case of $\mathscr{A}^M$ is the same.

From the concavity of the entropy function (Lemma \ref{lemma: subadditivity and concavity of entropy} (2)), 
for $\mu_n = \frac{1}{N}\sum_{n=0}^{N-1}T^n_*\sigma_N$
\begin{equation} \label{eq: estimate by concavity}
  H_{\mu_N}\left(\mathscr{B}^M\right) \geq \frac{1}{N}\sum_{n=0}^{N-1}H_{T^n_*\sigma_N}\left(\mathscr{B}^M\right)
     = \frac{1}{N}\sum_{n=0}^{N-1}H_{\sigma_N}\left(T^{-n}\mathscr{B}^M\right). 
\end{equation}     

Let $N = qM + r$ with $0\leq r <M$.
\begin{equation}  \label{eq: some edge effects}
   \begin{split}
   \sum_{n=0}^{N-1} H_{\sigma_N} \left(T^{-n}\mathscr{B}^M\right) & 
   = \sum_{t=0}^{M-1}\sum_{s=0}^q H_{\sigma_N}\left(T^{-sM-t}\mathscr{B}^M\right)
     - \sum_{n=qM+r}^{qM+M-1} H_{\sigma_N}\left(T^{-n}\mathscr{B}^M\right) \\
   & \geq    \sum_{t=0}^{M-1}\sum_{s=0}^q H_{\sigma_N}\left(T^{-sM-t}\mathscr{B}^M\right) - M \log \left|\mathscr{B}^M\right| \\
   & \geq    \sum_{t=0}^{M-1}\sum_{s=0}^q H_{\sigma_N}\left(T^{-sM-t}\mathscr{B}^M\right) - M^2 \log \left|\mathscr{B}\right|.
   \end{split}
\end{equation}   

We estimate $\sum_{s=0}^q H_{\sigma_N}\left(T^{-sM-t}\mathscr{B}^M\right)$ from below for each $t$.
We have 
\[ T^{-sM-t}\mathscr{B}^M = \bigvee_{m=0}^{M-1} T^{-(sM+t+m)}\mathscr{B}. \]
When we fix $0\leq t\leq M-1$ and move $0\leq s\leq q$ and $0\leq m \leq M-1$, the number $sM+t+m$ moves over 
\[  t, t+1, t+2, \dots, t+(q+1)M-1 \quad \text{without multiplicity}. \]
Hence 
\begin{equation*}
  \begin{split}
     \sum_{s=0}^q H_{\sigma_N}\left(T^{-sM-t}\mathscr{B}^M\right) &+ \sum_{n=0}^{t-1} H_{\sigma_N}\left(T^{-n}\mathscr{B}\right) \\
    & \geq   H_{\sigma_N}\left(\bigvee_{n=0}^{t+(q+1)M-1} T^{-n}\mathscr{B}\right)  \quad 
    \text{by Lemma \ref{lemma: subadditivity and concavity of entropy} (1)} \\
    &\geq H_{\sigma_N}\left(\mathscr{B}^N\right) \quad \text{by $t+(q+1)M \geq (q+1)M > N$}.
  \end{split}
\end{equation*}  
Therefore 
\begin{equation*}
  \begin{split}
     \sum_{s=0}^q H_{\sigma_N}\left(T^{-sM-t}\mathscr{B}^M\right) & \geq 
      H_{\sigma_N}\left(\mathscr{B}^N\right) - \sum_{n=0}^{t-1} H_{\sigma_N}\left(T^{-n}\mathscr{B}\right)  \\
     & \geq   H_{\sigma_N}\left(\mathscr{B}^N\right) - t \log \left|\mathscr{B}\right|  \\
     & \geq    H_{\sigma_N}\left(\mathscr{B}^N\right) - M \log \left|\mathscr{B}\right|  \quad 
     \text{by $t<M$}.
  \end{split}
\end{equation*}  
Thus 
\[ \sum_{t=0}^{M-1}\sum_{s=0}^q H_{\sigma_N}\left(T^{-sM-t}\mathscr{B}^M\right) \geq 
     M\cdot  H_{\sigma_N}\left(\mathscr{B}^N\right) - M^2 \log \left|\mathscr{B}\right|. \]
So by (\ref{eq: some edge effects})
\begin{equation*}
  \begin{split}
   \sum_{n=0}^{N-1} H_{\sigma_N}\left(T^{-n}\mathscr{B}^M\right) & \geq  
   \sum_{t=0}^{M-1}\sum_{s=0}^q H_{\sigma_N}\left(T^{-sM-t}\mathscr{B}^M\right) - M^2 \log \left|\mathscr{B}\right| \\
    & \geq    M\cdot  H_{\sigma_N}\left(\mathscr{B}^N\right) - 2 M^2 \log \left|\mathscr{B}\right|. 
  \end{split}
\end{equation*}   
From (\ref{eq: estimate by concavity}) we conclude that 
\[ \frac{1}{M}H_{\mu_N}\left(\mathscr{B}^M\right) \geq \frac{1}{N M}\sum_{n=0}^{N-1} H_{\sigma_N}\left(T^{-n}\mathscr{B}^M\right)
    \geq \frac{1}{N} H_{\sigma_N}\left(\mathscr{B}^N\right) - \frac{2M\log \left|\mathscr{B}\right|}{N}. \]
\end{proof}

We have 
\[  \int_X f\, d\mu_N = \frac{1}{N} \int_X \sum_{n=0}^{N-1} f\circ T^n \, d\sigma_N 
     = \frac{1}{N} \int_X \mathbb{S}_N f\, d\sigma_N.   \]
Claim \ref{claim: entropy bounded from below} implies 
\begin{equation*}
  \begin{split}
   & \frac{w}{M} H_{\mu_N}\left(\mathscr{B}^M\right) + \frac{1-w}{M}H_{\pi_* \mu_N}\left(\mathscr{A}^M\right) + 
    w \int_X f\, d\mu_N \\
    & \geq 
    \frac{w}{N} H_{\sigma_N}\left(\mathscr{B}^N\right) + \frac{1-w}{N} H_{\pi_*\sigma_N}\left(\mathscr{A}^N\right) 
    + \frac{w}{N} \int_X \mathbb{S}_N f \, d\sigma_N
    - \frac{2M\left(\log |\mathscr{A}| +\log |\mathscr{B}|\right)}{N} \\
    & = \frac{\log Z_N}{N} - \frac{2M\left(\log |\mathscr{A}| +\log |\mathscr{B}|\right)}{N}  \quad 
    \text{by Claim \ref{claim: Z_n and measure sigma_N}}.
    \end{split}
\end{equation*}    
Since $\mu_{N_k}\to \mu$ as $k\to \infty$, letting $N=N_k\to \infty$
\[ \frac{w}{M} H_{\mu}\left(\mathscr{B}^M\right) + \frac{1-w}{M}H_{\pi_* \mu}\left(\mathscr{A}^M\right) + 
    w \int_X f\, d\mu 
    \geq \lim_{N\to \infty} \frac{\log Z_N}{N}. \]
Here we have used the clopenness of the elements of $\mathscr{A}^M$ and $\mathscr{B}^M$.
Finally, letting $M\to \infty$, we get 
\[ w h_\mu\left(T, \mathscr{B}\right) + (1-w) h_{\pi_*\mu}\left(S, \mathscr{A}\right) 
    + w \int_X f\, d\mu \geq 
    \lim_{N\to \infty} \frac{\log Z_N}{N}. \]
\end{proof}

Now we can prove the main result (Theorem \ref{theorem: main theorem}).
We repeat the statement for the convenience of readers.

\begin{theorem}[$=$ Theorem \ref{theorem: main theorem}]
Let $\pi:(X, T)\to (Y, S)$ be a factor map between dynamical systems.
Then for any $0\leq w\leq 1$ and a continuous function $f:X\to \mathbb{R}$
\[  P^w(\pi, T, f) = P^w_{\mathrm{var}}(\pi, T, f). \]
\end{theorem}

\begin{proof}
We already proved in Proposition \ref{prop: topology bounds measure} that
\[ P^w_{\mathrm{var}}(\pi, T, f)  \leq P^w(\pi,T,f). \] 
By Corollary \ref{cor: reduction to zero dimensional case}, 
there exists a factor map $\pi^\prime:(X^\prime, T^\prime)\to (Y^\prime, S^\prime)$ between 
\textit{zero dimensional} dynamical systems with a continuous function 
$f^\prime:X^\prime\to \mathbb{R}$ such that 
\[ P^w(\pi,T, f) \leq P^w(\pi^\prime, T^\prime, f^\prime), \quad 
       P^w_{\mathrm{var}}(\pi^\prime, T^\prime, f^\prime) \leq P^w_{\mathrm{var}}(\pi,T,  f). \]
By Proposition \ref{prop: main result in zero dimensional case}    
\[ P^w(\pi^\prime, T^\prime, f^\prime) \leq P^w_{\mathrm{var}}(\pi^\prime, T^\prime, f^\prime). \]
Therefore 
\[        P^w(\pi,T, f) \leq P^w_{\mathrm{var}}(\pi,T,  f). \]
So we conclude that        
\[   P^w(\pi, T, f) = P^w_{\mathrm{var}}(\pi, T, f). \]
\end{proof}

\begin{remark}
The book of Downarowicz \cite{Downarowicz} systematically develops the idea of using 
zero dimensional dynamical systems in the study of entropy theory.
The above proof is influenced by this idea.
We also notice that it seems difficult to use this \textit{zero dimensional trick} in the proof of 
Proposition \ref{prop: topology bounds measure} in \S \ref{section: topology bounds measure}
because it is difficult to prove that principal extensions preserve weighted topological pressure without 
using the variational principle. A similar remark is given in \cite[Remark 7.6.12]{Downarowicz} about the proof of 
the standard variational principle.
\end{remark}

\vspace{0.5cm}

\address{ Masaki Tsukamoto \endgraf
Department of Mathematics, Kyushu University, Moto-oka 744, Nishi-ku, Fukuoka 819-0395, Japan}

\textit{E-mail}: \texttt{tsukamoto@math.kyushu-u.ac.jp}

\end{document}